\documentclass[a4paper]{article}
\pdfoutput=1 
\usepackage{hyperref}
\usepackage[english]{babel}
\usepackage[utf8]{inputenc}
\usepackage{amsmath}
\usepackage{amssymb}
\usepackage{amsthm}
\usepackage{graphicx}
\usepackage[colorinlistoftodos]{todonotes}

\theoremstyle{plain}
\newtheorem{defi}{Definition}[section]
\newtheorem{lemme}[defi]{Lemma}
\newtheorem{coro}[defi]{Corollary}
\newtheorem{prop}[defi]{Proposition}
\newtheorem{theo}[defi]{Theorem}
\theoremstyle{definition}
\newtheorem{example}[defi]{Example}
\theoremstyle{remark}
\newtheorem{remark}[defi]{Remark}
\numberwithin{equation}{section}

\def\nn{\mathbb{N}}
\def\zz{\mathbb{Z}}

\def\cc{\mathbb{C}}

\def\M{\mathfrak{M}}
\def\N{\mathfrak{N}}
\def\tM{\widetilde{\M}}
\def\L{\mathfrak L}
\def\P{\mathcal P}
\def\tP{\widetilde{\P}}
\def\H{\mathcal H}
\def\h{\mathfrak h}
\def\pre{\noindent{\bf Proof:}}
\def\fin{$\Box$}
\def\tr{\mathrm{Tr}}
\def\id{\mathrm{Id}}
\def\e{\mathrm{e}}
\def\i{\mathrm{i}}
\makeatletter
\newcommand{\oset}[2]{%
  {\mathop{#2}\limits^{\vbox to -.5\ex@{\kern-\tw@\ex@
   \hbox{\scriptsize #1}\vss}}}}
\makeatother
\def\submod{\oset{{\tiny d}}{-}}
\def\addmod{\oset{{\tiny d}}{+}}
\def\submodn{\oset{{\tiny n}}{-}}
\def\addmodn{\oset{{\tiny n}}{+}}\def\ind{1\hspace{-0.27em}\mathrm{l}}
\def\pp{\mathbb{P}}
\def\ee{\mathbb{E}}
\def\rhoinv{\rho^{\mathrm{inv}}}
\def\eps{\varepsilon}

\def\vec#1{|#1\rangle}
\def\connect{\overset{\mathfrak M}{\rightarrow}}

\def\biconnect{\overset{\mathfrak M}{\leftrightarrow}}

\def\enc{\mathrm{Enc}}

\def\braket#1#2{\langle #1  , #2\rangle}
\def\ketbra#1#2{|#1\rangle\langle#2|}

\title{Open Quantum Random Walks:\\ reducibility, period, ergodic properties}
\author{Raffaella Carbone\\Dipartimento di Matematica dell’Universit\'a di Pavia\\ via Ferrata, 1, 27100 Pavia, Italy\\\texttt{raffaella.carbone@unipv.it}\\ \\Yan Pautrat\\Laboratoire de Math\'ematiques\\
Universit\'e Paris-Sud\\
91405 Orsay Cedex, France\\\texttt{yan.pautrat@math.u-psud.fr}}
\begin{document}
\bibliographystyle{abbrv}
\maketitle

{\noindent\bf Abstract.} We study the analogues of irreducibility, period, and communicating classes for open quantum random walks, as defined in \cite{APSS}. We recover results similar to the standard ones for Markov chains, in terms of ergodic behavior, decomposition into irreducible subsystems, and characterization of stationary states.
\medskip

\medskip

\noindent{\bf AMSSC:} 81P16, 46L55, 60J05, 47A35.

\section{Introduction}

Open quantum random walks were recently defined by Attal \textit{et al.} in \cite{APSS}. These processes have a simple definition, implementing a Markovian dynamics influenced by internal degrees of freedom, and can be useful to model a variety of phenomena: quantum algorithms (see \cite{SP}), transfer in biological systems (see~\cite{MSKPE}) and possibly quantum exclusion processes. In addition, a continuous-time version can be defined (see \cite{PelOQRW}). Therefore, open quantum random walks seem to be good quantum analogues of Markov chains.

The usefulness of (classical) Markov chains, however, comes not only from the vast number of situations they can model, but also from the many properties implied by their simple definition. A textbook description of Markov chains, for instance, can start with the notion of irreducibility, which is easily characterized through the connectedness of the associated graph, and implies mean-ergodic convergence in law if an invariant probability exists (which is the case when the state space is finite). The next notion, the aperiodicity of an irreducible chain, is not as easy to characterize, but has simple sufficient conditions (\textit{e.g.} the existence of loops) and implies convergence in law, at least when the state space is finite. Last, the notion of connected subsets of the initial graph allows one to decompose a Markov chain into irreducible ones, to characterize its invariant states as convex combinations of invariant states for restricted chains.

On the other hand, the only general properties of open quantum random walks proven so far are the central limit theorem for the position process (see~\cite{AGS}) and the general K\"ummerer-Maassen theorem for quantum trajectories (see~\cite{KuMa}). In the present paper we discuss an analogue of the above textbook description of Markov chains, for open quantum random walks. The non-commutative nature of the objects under study, and specifically the fact that the transition probabilities are replaced by operators acting on a Hilbert space, are the cause of higher mathematical complexity. Some intuitive aspects of classical Markov chains, however, fruitfully remain, and we can recover a vision of irreducibility, period, and accessibility, in terms of paths. This is of interest for the study of more general quantum Markov processes, as it gives indications on the relevant extensions of classical concepts, and on techniques of proofs of associated results. We view this as an additional justification for the study of open quantum random walks.

Our theory will be constructed starting from pre-existing tools:
\begin{itemize}
\item a notion of irreducibility for general positive maps on non-commutative algebras, together with an associated Perron-Frobenius theorem, that was developed by various authors in the late seventies and early eighties (\cite{AHK}, \cite{EHK}, \cite{WE}, \cite{Gro});
\item a notion of period, together with associated results on the peripheral spectrum, that were defined in the same setting by Groh (\cite{Gro}) and extended by Fagnola and Pellicer (\cite{FP});
\item some old and new inspiring ergodic results (see \cite{FV} and \cite{KuMa}) and a decomposition of the support of invariant states proposed more recently by Baumgartner and Narnhofer (\cite{BN}) for quantum discrete time processes acting on finite dimensional spaces. 
\end{itemize}
We briefly describe the structure of the article and the main contents. Section~\ref{section_OQRWs} recalls the definitions, notations and basic results regarding open quantum random walks from \cite{APSS}. We describe the two types of (classical) processes associated to an OQRW: the process ``with (repeated) measurement", commonly called ``quantum trajectory", and the process ``without measurement". Sections \ref{section_irreducibility} and~\ref{section_aperiodicity} discuss, respectively, irreducibility and aperiodicity for OQRWs. Both follow the same structure: they start by recalling standard definitions and properties of irreducibility or aperiodicity for positive maps on operator algebras; then study the application to the special case of OQRWs.
Some immediate consequences on the ergodic behavior of the evolution are underlined. Section~\ref{section_recurrence} applies the results of the previous two sections to obtain convergence properties for irreducible, or irreducible aperiodic, open quantum random walks, for both processes described in section \ref{section_OQRWs}, \textit{i.e.} ``with measurement"
and ``without measurement". Section \ref{section_nonirreducible} expands on reducible open quantum random walks, characterizing in different ways their irreducible components. The resulting decomposition can be seen as related to a ``quantum communication relation'' among vectors of the underlying Hilbert space. Section \ref{section_stationary} states the general form of stationary states for reducible open quantum random walks. Its central point is the full exploitation of some results from \cite{BN}, which we state and prove in full detail. Section \ref{section_extension} mentions a natural extension of open quantum random walks, which are strongly related to the quantum Markov chains defined by Gudder in \cite{Gudder}. For this extension we discuss without proof a characterization of irreducibility, periodicity, communication classes, and their consequences: as we will see, all previous results will remain with paths on a graph replaced by paths on a multigraph. We conclude with section \ref{section_examples}, which is dedicated to examples and applications. We start a study of translation-invariant open quantum random walks on $\zz^d$ continued in \cite{CP2}, and extending that of \cite{AGS}. We study examples which illustrate our most practical convergence results, namely Corollaries \ref{coro_rec1}, \ref{coro_rec2}, and \ref{coro_rec3}, as well as our decomposition result, Theorem \ref{theo_invariantstates}.

\paragraph{Acknowledgements} The authors wish to thank Stéphane Attal for providing perpetual impetus to this project, Matteo Gregoratti for the organization of a meeting in Milano that played an important role in the development of this article, and Clément Pellegrini for many enthusiastic discussions. RC also gratefully acknowledges the support of PRIN project 2010MXMAJR 
and GNAMPA project ``Semigruppi markoviani su algebre non commutative'', and YP the support of ANR project “HAM-MARK”, n${}^\circ$ANR-09-BLAN-0098-01.

\section{Open quantum random walks}\label{section_OQRWs}

In this section we recall basic results and notations about open quantum random walks. For a more detailed exposition of OQRWs and related notions we refer the reader to \cite{APSS}.

We consider a Hilbert space $\mathcal H$ of the form 
$\mathcal H = \bigoplus_{i\in V}\mathfrak h_i$ where $V$ is a countable set of vertices, and each $\mathfrak h_i$ is a separable Hilbert space (making $\mathcal H$ separable).
This is a generalization with respect to standard OQRWs where the space $\mathcal H$ is $\mathfrak h \otimes\ell^2(V)$, or equivalently $\mathfrak h_i = \mathfrak h$ for all $i\in V$. This generalization will be useful when we consider decompositions of OQRWs, especially in section~\ref{section_nonirreducible}.
We view $\mathcal H$ as describing the degrees of freedom of a particle constrained to move on $V$: the ``$V$-component" describes the spatial degrees of freedom (the position of the particle) while $\mathfrak h_i$ describes the internal degrees of freedom of the particle, when it is located at site $i\in V$. 

For clarity, whenever a vector $x\in\mathcal H$ belongs to the subspace $\mathfrak h_i$, we will denote it by $x\otimes\vec i$, and drop the (implicit) assumption that $x\in\mathfrak h_i$. Similarly, when an operator $A$ on $\mathcal H$ satisfies $\mathrm{Ker}\, A \subset \mathfrak h_j^\perp$ and $\mathrm{Ran}\, A\subset \mathfrak h_i$, we denote it by $A=L_{i,j}\otimes \ketbra ij$ where $L_{i,j}$ is viewed as an operator from $\mathfrak h_j$ to $\mathfrak h_i$. Therefore, for $i,j,k$ in $V$, we have the relation
\[ 
\big(L_{i,j}\otimes\ketbra ij\big) \ \big(x\otimes\vec k\big) = 
\left\{ \begin{array}{ll} 0 & \mbox{ if }j\neq k, \\ 
L_{i,j}\,x\otimes \vec i & \mbox{ if }j=k. 
\end{array}\right. \]
All of these notations are consistent with the special case of $\mathcal H = \mathfrak h \otimes \ell^2(V)$, and with the interpretation of $\mathcal H$ described above.
\smallskip

We consider a map on the space $\mathcal I _1(\mathcal H)$ of trace-class operators on $\mathcal H$,
\begin{equation}\label{eq_OQRW}
\mathfrak M \, : \, \rho \mapsto \sum_{i,j\in V}A_{i,j}\,\rho \,A_{i,j}^*
\end{equation}
where, for any $i,j$ in $V$, the operator $A_{i,j}$ is of the form $L_{i,j}\otimes |i\rangle \langle j|$ and the operators $L_{i,j}$ satisfy
\begin{equation}\label{eq_stochastic}
\forall j\in V\quad \sum_{i\in V}L_{i,j}^* L_{i,j}=\id,
\end{equation}
where the series is meant in the strong convergence sense. The $L_{i,j}$ are thought of as encoding both the probability of a transition from site $j$ to site $i$, and the effect of that transition on the internal degrees of freedom. Equation \eqref{eq_stochastic} therefore encodes the ``stochasticity" of the transitions $L_{i,j}$.

Clearly \eqref{eq_OQRW} defines a trace-preserving (TP) map $\mathcal I_1(\mathcal H)\to \mathcal I_1(\mathcal H)$, which is completely positive (CP), \textit{i.e.} for any $n$ in $\nn^*$, the extension $\M\otimes \id$ to $\mathcal I_1(\mathcal H)\otimes \mathcal B(\cc^n)$ is positive. In particular, such a map transforms states (\textit{i.e.} positive elements of $\mathcal I_1(\mathcal H)$ with trace one) into states. A completely-positive, trace-preserving map will be called a CP-TP map. We shall call a map $\M$ as defined by \eqref{eq_OQRW} an open quantum random walk, or OQRW.  Note that \eqref{eq_stochastic} implies that $\|\M\|=1$ as an operator on $\mathcal I_1(\mathcal H)$ (see Remark \ref{remark_TPnormone} below).

\begin{remark}\label{remark_extensions}
In our interpretation of $L_{i,j}$ above, it would be more precise to say that the transition from site $j$ to site $i$ is encoded by the completely positive map $\rho_j\mapsto L_{i,j}\, \rho_j \, L_{i,j}^*$. A natural extension would be to replace this with a more general completely positive map $\rho_j\mapsto \Phi_{i,j}(\rho_j)$. We recover the ``transition operation matrices" introduced by Gudder in \cite{Gudder}. This will be discussed in section~\ref{section_extension}.
\end{remark}

Let us recall that the topological dual $\mathcal I_1(\mathcal H)^*$ can be identified with $\mathcal B(\mathcal H)$ through the duality
\[(\rho,X)\mapsto \tr(\rho\, X).\]
\begin{remark}\label{remark_TPnormone}
Trace-preservation of a map $\Phi$ is equivalent to $\Phi^*(\id)=\id$. The adjoint $\Phi^*$ is then a positive, unital (\textit{i.e.} $\Phi^*(\id)=\id$) map on $\mathcal B (\mathcal H)$, and by the Russo-Dye theorem (\cite{RD}) one has $\|\Phi^*\|=\|\Phi^*(\id)\|$ so that $\|\Phi\|=\|\Phi^*\|=1$.
\end{remark}

\begin{defi}
We say that an open quantum random walk $\mathfrak M$ is finite if $V$ is finite and every $\mathfrak h_i$ is finite-dimensional.
\end{defi}

\begin{remark}\label{remark_finite}
If an open quantum random walk is finite, then $\mathfrak M^*(\id)=\id$ implies that $1$ is an eigenvalue of $\mathfrak M$. Since $\M$ preserves the trace and the positivity, this implies that there exists an invariant state.
\end{remark}

\begin{remark} As noted in \cite{APSS}, classical Markov chains can be written as open quantum random walks. More precisely, if the transition matrix is $P=(P_{i,j})$ then, taking $L_{i,j}=\sqrt{P_{j,i}}\, U_{i,j}$ with any $U_{i,j}$ satisfying $U_{i,j}^* U_{i,j}=\id_{\mathfrak h_j}$, will preserve states of the form $\sum_{i\in V} p_i \otimes |i\rangle\langle i|$, and the induced dynamics on the family~$(p_i)_{i\in V}$ will be described by the transition matrix $P$. However, if $\mathrm{dim}\,\mathfrak h_i>~1$ we will run into possible non-uniqueness problems \textit{e.g.} for the invariant states of $\M$ (see section \ref{section_nonirreducible}). We feel this is an artificial degeneracy, not related to the properties of the Markov chain, but rather to the choice of the dilation. We will therefore only consider \textit{minimal} dilations of classical Markov chains, where $\mathrm{dim}\,\mathfrak h_i =1$ for all $i$ in $V$, and $L_{i,j}=\sqrt{P_{j,i}\,}$.
\end{remark}

A crucial remark is that, for any initial state $\rho$ on $\mathcal H$, which can be expanded~as
\[ \rho=\sum_{i,j\in V} \rho(i,j)\otimes \ketbra ij\]
and, for any $n\geq 1$, the evolved state $\mathfrak M ^n(\rho)$ is of the form
\begin{equation}\label{eq_Mn1}
\mathfrak M ^n (\rho) = \sum_{i\in V} \mathfrak M^n (\rho,i) \otimes \ketbra ii,
\end{equation}
where \textit{e.g.} for $n=1$,
\begin{equation}\label{eq_Mn2}
\mathfrak M^1 (\rho,i)=\sum_{j\in V} L_{i,j}\, \rho(j,j)\, L_{i,j}^*.
\end{equation}
Each $\mathfrak M^n(\rho,i)$ is a positive, trace-class operator on $\mathfrak h _i$ and
$ \sum_{i\in V} \tr\,\mathfrak M^n(\rho,i) =1$.
Therefore, the range of $\M$ is included in the set  ${\cal I_D}$ of block diagonal trace-class operators,
$$
{\cal I_D}=\{ \rho=\sum_{i\in V}\rho(i)\otimes \ketbra ii, \, \sum_{i\in V}\tr(|\rho(i)|)<+\infty  \},
$$
and $\mathcal I_D^*$ can be identified with 
\[\mathcal B_D =\{ X=\sum_{i\in V} X(i)\otimes \ketbra ii, \ \mathrm{sup}\,\|X(j)\|_{\mathcal B(\mathfrak h_j)}<\infty \}.\]
This feature will have a great importance in the characterization of many properties of OQRWs, \textit{e.g.}:
\begin{enumerate}
\item the invariant states of an OQRW $\mathfrak M$ belong to $\cal I_D$,
\item the reducibility of $\mathfrak M$ can be established considering only block-diagonal projections (see section \ref{section_irreducibility}) and the only meaningful enclosures are generated by vectors of the form $x\otimes \vec i$  (see section \ref{section_nonirreducible}),
\item the cyclic projections defining the period have block-diagonal form (see section~\ref{section_aperiodicity}).
\end{enumerate}
In addition, we remark from  \eqref{eq_Mn2}  that $\M^n(\rho)$ depends only on the diagonal elements $\rho(i,i)$. Therefore, from now on, we will only consider states of the form $\rho=\sum_{i\in V}\rho(i)\otimes \ketbra ii$. Equation \eqref{eq_Mn2} remains valid, replacing $\rho(i,i)$ by~$\rho(i)$.
\smallskip

We now describe the (classical) processes of interest, associated with $\M$. We start from a state $\rho$ which we assume to be of the form $\rho=\sum_{i\in V}\rho(i) \otimes \ketbra ii$. We evolve $\rho$ for a time $n$, obtaining the state $\M^n(\rho)$ as in \eqref{eq_Mn1}. We then make a measurement of the position observable. According to standard rules of quantum measurement, we obtain the result $i\in V$ with probability $\tr\,\M^n(\rho,i)$. Therefore, the result of this measurement is a random variable $Q_n$, with law $\pp(Q_n=i)= \tr\,\mathfrak M^n(\rho,i)$ for $i\in V$. In addition, if the position $Q_n=i\in V$ is observed, then the state is transformed to $\frac{\mathfrak M^n(\rho,i)}{\tr\,\mathfrak M^n(\rho,i)}$. We call this process $(Q_n,\frac{\M^n(\rho,Q_n)}{\tr\,\M^n(\rho,Q_n)})_n$ the process ``without measurement" to emphasize the fact that virtually only one measurement is done, at time $n$. Notice that, in practice, two values of this process at times $n<n'$ cannot be considered simultaneously since the measure at time $n$ perturbs the system, and therefore subsequent measurements. This is reflected in the fact that \textit{a priori}, $\M^n(\rho)$ and $\M^{n'}(\rho)$ do not commute (see~\cite{DaviesBook} for a short introduction to the conceptual difficulties of associating random variables to operators).

Now assume that we make a measurement at every time $n\in\nn$, applying the evolution by $\M$ between two measurements. Again assume that we start from a state $\rho$ of the form $\sum_{i\in V}\rho(i)\otimes |i\rangle \langle i|$. Suppose that at time $n$, the position was measured at $X_n=j$ and the state (after the measurement) is $\rho_n \otimes \ketbra {j}{j}$. Then after the evolution, the state becomes
\[\M(\rho_n\otimes \ketbra jj)= \sum_{i\in V} L_{i,j}\, \rho_n\, L_{i,j}^* \otimes \ketbra ii\]
so that a measurement at time $n+1$ gives a position $X_{n+1}=i$ with probability $\tr\, L_{i,j}\, \rho_n\, L_{i,j}^*$, and then the state becomes $\rho_{n+1}\otimes\ketbra ii$ with $\rho_{n+1}=\frac{L_{i,j}\, \rho_n\, L_{i,j}^*}{\tr\, L_{i,j}\, \rho_n\, L_{i,j}^*}$. The sequence of random variables $(X_n,\rho_n)$ is therefore a Markov process with transitions defined by
\[ \pp\Big((X_{n+1},\rho_{n+1})\!=\!(i,\!\frac{L_{i,j}\,\rho_n\, L_{i,j}^*}{\tr(L_{i,j}\rho L_{i,j}^*)}\!)\Big|(X_n,\rho_n)\!=\!(j,\rho_n)\Big)=\tr (L_{i,j}\,\rho_n\, L_{i,j}^*)\quad \forall i\in V,\]
and initial law $\pp\big((X_0,\rho_0)=(i,\frac{\rho(i)}{\tr\rho(i)})\big)=\tr\,\rho(i).$
Note that the sequence of positions $X_0~=~i_0$, \ldots, $X_n=i_n$ is observed with probability 
\[\tr\,L_{i_n,i_{n-1}}\ldots L_{i_1,i_0}\,\rho(i_0)\, L_{i_1,i_0}^*\ldots L_{i_n,i_{n-1}}^*\]
and completely determines the state $\rho_n$:
\begin{equation}\label{eq_rhon} \rho_n = \frac{L_{i_n,i_{n-1}}\ldots L_{i_1,i_0}\,\rho(i_0)\, L_{i_1,i_0}^*\ldots L_{i_n,i_{n-1}}^*}{\tr\,L_{i_n,i_{n-1}}\ldots L_{i_1,i_0}\,\rho(i_0)\, L_{i_1,i_0}^*\ldots L_{i_n,i_{n-1}}^*}.
\end{equation}
As emphasized in \cite{APSS}, this implies that for every $n$ the laws of $X_n$ and $Q_n$ are the same, \textit{i.e.} \[ \pp(X_n=i)=\pp(Q_n=i)\quad \forall i\in V.\]
It also implies for any $n$ the relation 
\begin{equation}\label{eq_exprhon}
\ee(\rho_n\otimes\ketbra {X_n}{X_n})=\M^n(\rho).
\end{equation}

\section{Irreducibility for OQRWs}\label{section_irreducibility}
In this and in the following sections, $\Phi$ is assumed to be a positive map on the ideal $\mathcal I_1(\mathcal H)$ of trace operators on some given Hilbert space $\mathcal H$. We recall that such a map is automatically bounded as a linear map on $\mathcal I_1(\mathcal H)$ (see \textit{e.g.} Lemma~2.2 in \cite{Sch}), so that it is also weak-continuous. In most practical cases, we will additionally assume that $\|\Phi\|=1$; as we noted in Remark \ref{remark_TPnormone}, this will be the case, in particular, if~$\Phi$ is trace-preserving.

We recall some standard notations: an operator $X$ on $\mathcal H$ is called positive, denoted $X\geq0$, if for $\phi\in \mathcal H$,  one has $\langle \phi, X\, \phi\rangle \geq 0$. It is called definite positive, denoted $X>0$, if for $\phi\in \mathcal H\setminus\{0\}$, one has $\langle\phi, X\, \phi\rangle>0$.
\smallskip

We summarize here the definition of irreducibility introduced by Davies (see~\cite{Dav}), and other related notions. We shall see in Proposition \ref{ergodic=irreducible} that the first two (irreducibility and ergodicity) are equivalent.

\begin{defi}
The positive map $\Phi$ is called:
\begin{itemize}
\item irreducible if the only orthogonal projections $P$ reducing $\Phi$, \textrm{i.e.} such that~$\Phi\big(P\mathcal I_1(\mathcal H)P\big) \subset P\mathcal I_1(\mathcal H)P$, are $P=0$ and $\id$,
\item ergodic if, for any $\rho\geq 0$, $\rho\neq 0$ in $\mathcal I_1(\mathcal H)$, there exists $t$ such that~$\e^{t \Phi} (\rho)>~0$,
\item positivity-improving if, for any $\rho\geq 0$, $\rho\neq 0$ in $\mathcal I_1(\mathcal H)$, one has
$\Phi(\rho)>0$,
\item
$n$-regular for $n\in{\mathbb N}^*$ if $\Phi^n$ is positivity improving.
\end{itemize}
\end{defi}

\begin{remark}
The condition $\Phi\big(P\mathcal I_1(\mathcal H)P\big) \subset P\mathcal I_1(\mathcal H)P$ is equivalent to the condition $\Phi(P)\leq \alpha P$ for some $\alpha>0$ whenever $P\in \mathcal I_1(\mathcal H)$, i.e. whenever $P$ is finite-dimensional. In the infinite-dimensional case one can prove that $P$ reduces $\Phi$ if and only if for any finite-dimensional projection $Q$ with $Q\leq P$, one has~$\Phi(Q)\leq \alpha  P$.
\end{remark}

\begin{remark}
An equivalent formulation of ergodicity is that for any $\rho\geq 0$, $\rho\neq 0$ in $\mathcal I_1(\mathcal H)$, for any $t>0$ one has $\e^{t \Phi} (\rho)>0$. This follows from the observation that the support projection of $\e^{t \Phi}(\rho)$ does not depend on $t>0$. 
\end{remark}

There is a possible confusion here due to the fact that some authors (\cite{FP},
\cite{Gro}) work in the Heisenberg representation, \textit{i.e.} in our notation consider $\Phi^*$, while others (\cite{EHK}, \cite{Sch}), like us, work in the Schr\"odinger representation. For completeness we give the next proposition, which connects the two representations:

\begin{prop} 
Let $\Phi$ be a positive, trace-preserving map on~$\mathcal I_1(\mathcal H)$.
\begin{itemize}
\item An orthogonal projection $P$ reduces $\Phi$ if and only if $P\le \Phi^*(P)$, \textit{i.e.} $1-P$ reduces~$\Phi^*$,
\item $\Phi$ is ergodic (resp. positivity improving, regular) if and only if $\Phi^*$ is ergodic (resp. positivity improving, regular).
\end{itemize}
\end{prop}

\pre

Assume first that $P$ reduces $\Phi$, i.e.  
$\Phi(P{\cal I}_1({\cal H})P)\subset P{\cal I}_1({\cal H})P$. Then for any~$\phi\in\mathcal H$ of norm one,
\[\langle \phi,P \phi \rangle = \tr (\ketbra{P\phi}{P\phi}) =\tr\big(\Phi(|P\phi\rangle \langle P\phi|)\big)\] 
and, by the reduction assumption, this is
\[\tr\big(P\,\Phi(|P\phi\rangle \langle P\phi|)\big) = \langle P\phi,\Phi^*(P) P\phi \rangle \leq   \langle \phi,\Phi^*(P) \, \phi \rangle,\]
so that $P\leq \Phi^*(P)$. Conversely, if $P\le\Phi^*(P)$, then, for any trace-class $\rho\geq 0$,
\[
\tr(P\rho P) \leq \tr(P\rho P\,\Phi^*(P))
= \tr(P\,\Phi(P\rho P)\,P)
\le \tr(\Phi(P\rho P)) = \tr(P\rho P).
\]
We therefore have the equality
$\tr(P\,\Phi(P\rho P)\,P)
= \tr(\Phi(P\rho P))$  which implies the inclusion $\Phi(P\rho P)\in P{\cal I}_1({\cal H})P$ for $\rho\geq 0$, hence for any $\rho \in\mathcal I_1(\mathcal H)$.

To prove the second point consider \textit{e.g.} $\Phi$ ergodic. For any $\rho\geq 0$, $\rho\neq 0$ in~$\mathcal I_1$ one has
$\e^{t \Phi} (\rho)>0$ for all $t>0$. So, for any bounded positive, non-zero operator~$X$, we have
$
0 < \tr(X\e^{t \Phi} (\rho)) = \tr(\e^{t \Phi^*} \!(X)\,\rho).
$
Taking $\rho$ of the form~$|\phi\rangle\langle \phi|$, we deduce $\e^{t \Phi^*} (X)>0$.
Other statements are proved in the same way.
\fin

The article \cite{EHK} shows that irreducibility and ergodicity are equivalent, but considers only the finite-dimensional case. We extend this statement to the infinite-dimensional case below:
\begin{prop}\label{ergodic=irreducible}
A positive map $\Phi$ on $\mathcal I_1(\mathcal H)$ is ergodic if and only if it is irreducible.
\end{prop}
\pre

If $\Phi$ is not irreducible, then there exists a non-trivial projection $P$ and a non-negative trace-class operator $\rho$ such
that $\Phi(P\rho P)\leq \alpha P$ and $P\rho P\neq 0$ but then, for any $t$, one has $\e^{t\Phi}(P\rho P)\leq \e^{t\alpha} P$
so that $\e^{t\Phi}(P\rho P)$ is non-definite for all $t$ and $\Phi$ is not ergodic.

To prove the converse we use the characterization in terms of the dual $\Phi^*$. Assume $\Phi$, hence $\Phi^*$, is irreducible, consider $X\geq 0$, $X\neq 0$ in $\mathcal B(\mathcal H)$; for a fixed~$t>0$ let
$$ e_p(X)= \sum_{k=0}^p \frac{t^k}{k!} \Phi^{*k}(X).$$
Define $P$ to be the support projection of $\e^{t\Phi^*}(X)$ and $P_{p}=\ind_{[1/p,+\infty[}(e_p(X))$.
Obviously $P_p\leq P$ and $P_p \leq p \, e_p(X)$ for all $p$, and $P_p \to P$ in the sense of strong convergence as $p\to\infty$, thanks to the properties of bounded measurable functional calculus (see \textit{e.g.} Theorem VII.2 in \cite{RS1}). We have:
\begin{eqnarray*} \frac1p \,\Phi^*(P_{p}) & \leq & \Phi^*(e_p(X)) 
= \sum_{k=1}^{p+1} \frac kt \frac{t^k}{k!} \,\Phi^{*k}(X)\\
&\leq& \frac{p+1}t \,e_{p+1}(X) \leq \frac{p+1}t \, \e^{t\Phi^*}(X)
\end{eqnarray*}
so that $\mathrm{supp}\,\Phi^*(P_{p})\subset \mathrm{supp}\,P$ and, by the weak-$\ast$ continuity of $\Phi^*$, one has $\mathrm{supp}\,\Phi^*(P)\subset\mathrm{supp}\, P$, \textit{i.e.} $P$ reduces $\Phi^*$. Since $\e^{t\Phi^*}(X)\geq X$, the projector $P$ cannot be zero, so by irreducibility $P$ is $\id$ and $\e^{t\Phi^*}(X)>0$.
\fin

\begin{remark} In \cite{EHK}, ergodicity is defined for finite dimensional $\mathcal H$
by the property that $(\id+\Phi)^{\mathrm{dim}\,\mathcal{H}-1}$ is positivity-improving. This is equivalent to our
definition: see the remark following Lemma 3.1 in \cite{Sch}.
\end{remark}
\smallskip

When speaking about reducibility/irreducibility of quantum maps, one enters a jungle of different approaches and terminologies, which, in many cases, are essentially equivalent.
Concerning this, we recall that a reducing projection~$P$ is called by some authors a \textit{subharmonic projection} for $\Phi^*$, following the line common to the classical literature on Markov chains.

Also, more recently (in \cite{BN}, as far as we know), the notion of enclosure has been introduced in the context of CP-TP maps.
A closed vector space $\mathcal V$ is called an \textit{enclosure} if~$\mathrm{supp} \, \rho\subset \mathcal V$ implies $\mathrm{supp} \, \mathfrak M(\rho)\subset \mathcal V$. It is immediate that a space $\mathcal V$ is an enclosure if and only if the projection $P$ on $\mathcal V$ reduces~$\mathfrak M$. So, an equivalent way to define irreducibility is asking that there exist no non-trivial enclosures. The notion of enclosure will be crucial in the discussion of decompositions of reducible open quantum random walks (see section \ref{section_nonirreducible}).

\medskip

Next, we characterize irreducibility in terms of unravellings. We consider a completely positive trace preserving map $\Phi$ and fix an unravelling $(A_\kappa)_{\kappa\in K}$ of~$\Phi$, provided by Kraus' representation theorem (see \cite{Kraus} or \cite{NieChu}, where this is called the operator-sum representation):
\begin{equation}\label{eq_Kraus}\Phi(\rho)=\sum_{\kappa\in K} A_\kappa \rho A_\kappa^\ast.\end{equation}
We will characterize  irreducibility (and the property of being positivity-improving or $n$-regular) in terms of an unravelling $(A_\kappa)_{\kappa\in K}$. We denote by $\cc[A]$ the set of polynomials in $A_\kappa$, \textit{i.e.}
the algebra (\textit{not} the *-algebra ) generated by the operators $A_\kappa$, $\kappa\in K$.
The following result summarizes Schrader's Lemmas 3.3 and 3.4 (\cite{Sch}):
\begin{lemme}\label{lemma_ergodicityCPTPmaps}
A completely positive map $\Phi$ of the form \eqref{eq_Kraus} is :
\begin{itemize}
\item positivity improving if and only if for any $\phi\in\mathcal H\setminus\{0\}$, the set 
$\{A_\kappa \,\phi, \, \kappa\in~K\}$
is total in~$\mathcal H$,
\item 
$n$-regular if and only if for any $\phi\in\mathcal H\setminus\{0\}$, the set
$\{A_{k_1}\!\ldots\! A_{k_n} \phi, \, k_1,\ldots,k_n\in~K\}$
is total in $\mathcal H$,
\item irreducible if and only if for any $\phi\in\mathcal H\setminus\{0\}$, the set
$ \cc[A] \, \phi$
is dense in~$\mathcal H$.
\end{itemize}
\end{lemme}
Lemmas 3.3 and 3.4 in \cite{Sch} are stated in terms of the operators $A_\kappa^*$. The connection with our statement comes from the following straightforward lemma:
\begin{lemme}\label{lemma_totality}
Consider a family $(A_\kappa)_{\kappa\in K}$ of operators on $\mathcal H$. Then the following are equivalent:
\begin{itemize}
\item for any $\phi\in\mathcal H\setminus\{0\}$, the set $\{A_\kappa\,\phi, \kappa\in K\}$ is total in $\mathcal H$,
\item for any $\phi$ and $\psi$ in $\mathcal H\setminus\{0\}$, there exists $k$ in $K$ such that $\langle \psi, A_\kappa \,\phi\rangle\neq 0$,
\item for any $\phi\in\mathcal H\setminus\{0\}$, the set $\{A_\kappa^*\, \phi, \kappa\in K\}$ is total in $\mathcal H$.
\end{itemize}\end{lemme}

Before we state our characterization of irreducibility for open quantum random walks, let us introduce some notation: for $i,j$ in $V$ we call a \textit{path from $i$ to~$j$} any finite sequence $i_0,\ldots,i_\ell$ in $V$ with $\ell\geq 1$, such that $i_0=i$ and $i_\ell=j$. Such a path is said to be \textit{of length $\ell$}. We denote by $\mathcal P (i,j)$ (resp. $\mathcal P_\ell (i,j)$) the set of paths from $i$ to $j$ of arbitrary length (resp. of length $\ell$). A path from $i$ to $i$ will be called a loop; by convention we consider the sequence $\{i\}$ as a loop (with length one), \textit{i.e.} an element of $\mathcal P(i,i)$. For $\pi=(i_0,\ldots,i_\ell)$ in $\mathcal P (i,j)$ we denote by $L_\pi$ the operator from $\mathfrak h_i$ to $\mathfrak h_j$:
$$ L_\pi=L_{i_\ell,i_{\ell-1}}\ldots L_{i_1,i_0}=L_{j,i_{\ell-1}}\ldots L_{i_1,i}.$$

We can now prove:
\begin{prop}\label{prop_ergodicityOQRW}
The CP-TP map $\M$ is irreducible if and only if, for every $i$ and $j$ in $V$, one of the following equivalent conditions holds:
\begin{itemize}
\item for any $x$ in $\mathfrak h_i\setminus\{0\}$, the set
$\{L_\pi x \, |\, \pi\in\mathcal P(i,j)\}$
is total in $\mathfrak h_j$,
\item for any $x$ in $\mathfrak h_i\setminus\{0\}$ and $y$ in $\mathfrak h_j\setminus\{0\}$ there exists a path $\pi$ in $\mathcal P(i,j)$ such that~$\langle y, L_\pi x\rangle~\neq~0.$
\end{itemize}
\end{prop}
\pre

This is an immediate application of Lemmas \ref{lemma_ergodicityCPTPmaps} and \ref{lemma_totality}, and the observation that, if $A_{j,i}=L_{j,i}\otimes\ketbra ji$, then
\[ A_{j_\ell,i_\ell}\ldots A_{j_1,i_1}=\left\{ \begin{array}{ll} L_{j_\ell,i_{\ell}}\ldots L_{i_2,i_1} \otimes \ketbra {j_\ell}{i_1}  & \mbox{ if }i_\ell=j_{\ell-1},\ldots, i_2=j_1,\\ 0 & \mbox{ otherwise.}\end{array}\right. \Box \]
\smallskip

As a first application we prove that ``positivity-improving" is an essentially useless notion in the framework of open quantum random walks, and that we have a constraint on the values $n$ allowing $n$-regularity:

\begin{coro}
The CP-TP map $\M$ is positivity-improving if and only if every~$\mathfrak h_i$ is one-dimensional and the underlying classical Markov process has positive transition probabilities.

It is $n$-regular if and only if, for any $i$ and $j$ in $V$, one of the equivalent formulations holds:
\begin{itemize}
\item for any nonzero $x$ in $\mathfrak h_i$, the set
$\{L_\pi x \, |\, \pi\in\mathcal P_n(i,j)\}$
is total in $\mathfrak h_j$,
\item for any nonzero $x$ in $\mathfrak h_i$ and $y$ in $\mathfrak h_j$ there exists a path $\pi$ in $\mathcal P_n(i,j)$ such that $\langle y, L_\pi x\rangle \neq 0.$
\end{itemize}
 A necessary condition for $n$-regularity is $\mathrm{card}\, \mathcal P_n(i,j)\geq \mathrm{dim}\, \mathfrak h_j$ for all $i,j\in V$. 
\end{coro}

\pre

By the first point of Lemma \ref{lemma_ergodicityCPTPmaps}, $\M$ is positivity-improving iff, for all $i,j\in V$, the set $\{L_{j,i}\,x\}$ is total in $\mathfrak h_j$ for any $x$ in $\mathfrak h_i$.
Since $\mathrm{dim}\,\mathrm{Vect}\{L_{j,i}\,x\} \leq 1$,
this is possible only if $\mathrm{dim}\, \mathfrak h_j=1$. In that case, the open quantum random walk is the minimal dilation of a classical Markov chain and the statement is obvious. The other statements are obtained by applying these requirements to $\mathfrak M^n$.
\fin

\smallskip
We can therefore give the following definition for an irreducible OQRW, which emphasizes our interpretation in terms of paths.
\begin{defi}\label{defi_irroqrw}
Let $\mathfrak M$ be an open quantum random walk. We say that two sites $i,j$ in $V$ are connected by $\mathfrak M$, which we denote by $i\connect j$, if one of the equivalent conditions of Proposition \ref{prop_ergodicityOQRW} holds. As we have shown, $\mathfrak M$ is irreducible if and only if, for any two $i$ and $j$ in $V$, one has $i\connect j$ and $j\connect i$.
\end{defi}

\begin{remark} A minimal dilation of a classical Markov chain is irreducible if and only if the Markov chain is irreducible in the classical sense.
\end{remark}

Until now, we have basically found necessary and sufficient conditions for irreducibility of an open quantum random walk. In section \ref{section_nonirreducible} we will discuss decompositions of reducible open quantum random walks into irreducible ones.
\smallskip

The following proposition essentially comes from \cite{Sch}:
\begin{prop}\label{prop_Schrader}
Assume a 2-positive map $\Phi$ on $\mathcal I_1(\mathcal H)$ has an eigenvalue $\lambda$ of modulus $\|\Phi\|$, 
with eigenvector $\rho$. Then:
\begin{itemize}
\item $\|\Phi\|$ is also an eigenvalue, with eigenvector $|\rho|$,
\item if $\Phi$ is irreducible, then $\lambda$ is a simple eigenvalue.
\end{itemize}
In particular, if $\Phi$ is irreducible and has an eigenvalue of modulus $\|\Phi\|$,
then~$\|\Phi\|$ is a simple eigenvalue, with an eigenvector that is definite-positive.
\end{prop}

\begin{remark}
Here and in the rest of this paper, by a simple eigenvalue of an operator $\Phi$ we mean a scalar $\lambda$ such that $\mathrm{dim\,Ker}\,(\Phi-\lambda\,\id)=1$.
\end{remark}

\pre

Theorems 4.1 and 4.2 from \cite{Sch} give us the first two statements. The third one follows from the fact that $\exp\|\Phi\|\times |\rho|=(\exp\Phi)(|\rho|)>0$ by irreducibility. \fin
\smallskip

In relation with the above results we can prove the following ergodic convergence
result for irreducible 2-positive, trace-preserving maps, which applies in particular to the case of CP-TP maps and can be seen as a discrete time version of the Frigerio-Verri ergodic theorem (\cite{FV}, Theorem 1.1):
\begin{prop}\label{prop_ergodicconvergence}
Let $\Phi$ be a positive contraction on $\mathcal I_1(\mathcal H)$ that has~$1$ as a simple eigenvalue. Then the associated eigenvector is  (up to normalization) an invariant state $\rhoinv$ and, for any state $\rho$, one has the weak convergence
\begin{equation}\label{eq_propergodicconvergence}
 \frac1n\sum_{k=0}^{n-1} \Phi^k(\rho)\underset{n\to\infty}{\longrightarrow} \rhoinv.
\end{equation}
\end{prop}

\pre

Consider an invariant trace-class operator $\rhoinv$. Since $\Phi$ preserves positivity, one can assume that $\rhoinv \geq 0$ and by necessity its trace is non-zero, so it can be assumed to have trace one. Define $\Psi_n=\frac1n\sum_{k=0}^{n-1}\Phi^k$. One has
$ \tr[\Psi_n(\rho)\,X]=\tr[\rho \,\Psi_n^*(X)]$. By the Banach-Alaoglu theorem, $\Psi_n^*(X)$ 
has weak-$\ast$ convergent subsequences. Denote by $Y$ the weak-$\ast$ limit of a 
subsequence $\Psi_{n_k}^*(X)$; one has $\Phi^*\circ\Psi_{n_k}^*(X)\to \Phi^*(Y)$ by the weak-$\ast$-continuity of $\Phi^*$, and, for any trace-class~$\rho$,
\begin{eqnarray*}
\tr[\rho(\id-\Phi^*)(Y)]&=&\lim_k\tr[\rho\,(\id-\Phi^*)\,\Psi_{n_k}^*(X)]\\
&=&
 \lim_k \tr[\rho\, \frac1{n_k}\sum_{j=0}^{n_k-1}(\Phi^{*j}-\Phi^{*(j+1)})(X)] 
\\
&=&\lim_k\tr[\rho\, \frac1{n_k}(\id-\Phi^{*{n_k}})(X)]  = 0,
\end{eqnarray*}
so that $\tr\left[(\id-\Phi)(\rho)\, Y\right]=0$, for any $\rho$, implying $Y\in\left(\mathrm{Ran}\, (\id-\Phi)\right)^\perp$. That space is of dimension one by assumption, so that $Y=\lambda_X \id$ and we have $\lim_k \tr\left[\Psi_{n_k}(\rho)\,X\right]\to \lambda_X$ for any trace-class $\rho$.  Writing this for $\rho$ equal to the eigenvector $\rhoinv$ leads to $\lambda_X=\tr(\rhoinv X)$, showing that $\lambda$ is independent on the subsequence $(n_k)_k$. When $\rho$ is a state we obtain the convergence \eqref{eq_propergodicconvergence}. This concludes the proof.
\fin

The following theorem is a direct application of Proposition \ref{prop_Schrader}:
\begin{theo}\label{theo_unicity}
An irreducible open quantum random walk $\M$ has an invariant state if and only if $1$ is an eigenvalue of $\mathfrak M$. If it does, then it has only one, and that invariant state is faithful.
\end{theo}

A second theorem follows from Proposition \ref{prop_ergodicconvergence}:
\begin{theo}\label{theo_ergodicconvergence}
Assume that an open quantum random walk $\M$ is irreducible and has an invariant state~$\rhoinv$. For any state $\rho$, one has
$\frac 1n\sum_{k=0}^{n-1}\mathfrak M ^k(\rho) \to \rhoinv$ weakly.
\end{theo}

\section{Period and aperiodicity for OQRWs}\label{section_aperiodicity}

As in the previous section, we start with a review of the notion of period for a positive trace-preserving map $\Phi$. Here we follow Fagnola and Pellicer (\cite{FP}) and Groh (\cite{Gro}). We define $\submod$ to be subtraction \textit{modulo} $d$.

\begin{defi}\label{def-period}
Let $\Phi$ be a positive, trace-preserving, irreducible map and let
$(P_0,\ldots,P_{d-1})$ be a resolution of identity, \emph{i.e.} a family of orthogonal
projections such that $\sum_{k=0}^{d-1} P_k=\id$. One says that $(P_0,\ldots,P_{d-1})$
is $\Phi$-cyclic if $\Phi^*(P_k)=P_{k\submod 1}$ for $k=0,\ldots,d-1$. 
The supremum of all $d$ for which there exists a $\Phi$-cyclic resolution of identity $(P_0,\ldots,P_{d-1})$ is called the \textrm{period} of $\Phi$.
If $\Phi$ has period~$1$ then we call it aperiodic.
\end{defi}

\begin{remark}
Even if in an embryonic stage, we recall that a characterization of a cyclic resolution of the identity was already given, in the Schr\"odinger picture, in \cite{EHK}, Theorem 3.4. 
\end{remark}

\begin{example}\label{ex_Orey}
Define a quantum Orey chain to be a CP-TP map $\Phi$ on $\mathcal I_1(\mathcal H)$ such that for any $\rho,\sigma$ in $\mathcal I_1(\mathcal H)$ one has $\Phi^n(\rho)-\Phi^n(\sigma)\rightarrow0$ in trace-norm, as~$n\to\infty$. A quantum Orey chain is aperiodic. Indeed, if $(P_0,\ldots,P_{d-1})$ is a cyclic resolution of identity with $d\geq 1$ then for $\rho$, $\sigma$ satisfying 
$ \mathrm{supp}\,\rho\subset \mathrm{Ran}\, P_0$, $\mathrm{supp}\,\sigma\subset \mathrm{Ran}\, P_1$, we have 

$$
\tr{\big((\Phi^n(\sigma)-\Phi^n(\rho))P_k\big)}
=\tr{\big((\sigma-\rho)\Phi^{*n}(P_k)\big)}
=\left\{\begin{array}{ll}
1 & \mbox{ for } k\submod n = 0 \\
-1 & \mbox{ for } k\submod n = 1 \\
0 & \mbox{ otherwise, } 
\end{array}\right.
$$ 
which contradicts the Orey property.
\end{example}

The following result is a combination of Theorems 3.7 and 4.3 of Fagnola-Pellicer in \cite{FP} (the latter was also partially proven by Groh in \cite{Gro}). Note that these results are proven in finite dimension, but they immediately extend to infinite dimension.
\begin{prop}\label{prop_aperiodicCPTP}
If $\Phi$ is an irreducible, 2-positive map on $\mathcal I_1(\mathcal H)$ and has finite period $d$ then:
\begin{itemize}
\item the peripheral point spectrum of $\Phi^*$, i.e. the set $\mathrm{Sp}_{pp}\Phi^* $, is a subgroup of the circle group $\mathbb T$,
\item the primitive root of unity $\e^{\i 2\pi/d}$ belongs to $\mathrm{Sp}_{pp}\Phi^*$ if and only if $\Phi$ is~$d$-periodic.
\end{itemize}
\end{prop}
The following result is an immediate consequence of Proposition \ref{prop_aperiodicCPTP}.
\begin{prop}\label{prop_cvgCPTPaperiodic}
If a $2$-positive TP map $\Phi$ on $\mathcal I_1(\mathcal H)$ is irreducible and aperiodic with invariant state $\rhoinv$, and $\mathcal H$ is finite-dimensional then
\begin{itemize}
\item $\mathrm{Sp_{pp}}\,\Phi\cap\mathbb T = \{1\}$ 
\item for any $\rho\in \mathcal I_1(\mathcal H)$ one has $\Phi^n(\rho)\rightarrow \rhoinv$ as $n\to\infty$.
\end{itemize}
\end{prop}

When considering completely positive maps, we will need to be able to characterize the period in terms of an unravelling
to apply it to OQRWs. We therefore fix an unravelling $A=(A_\kappa)_{\kappa\in K}$ of $\Phi$, \textit{i.e.} $\Phi(\rho)=\sum_{\kappa\in K}A_\kappa\,\rho\, A_\kappa^*$.
\begin{defi} 
Let $(P_0,\ldots,P_{d-1})$ be a resolution of identity. One says that it
is $A$-cyclic if $P_j A_\kappa= A_\kappa P_{j\submod 1}$ for $j=0,\ldots,d-1$ and any $k$.
\end{defi}

The following is Theorem 5.4 from \cite{FP}, which again extends to the infinite dimensional case, with the same proof.
\begin{prop}\label{prop_Vcyclic}
Let $\Phi$ be an irreducible CP-TP map on $\mathcal I_1(\mathcal H)$.
A resolution of the identity is $\Phi$-cyclic if and only if it is $A$-cyclic.
\end{prop}

\begin{remark}
For a $\Phi^*$-invariant weight (not necessarily a state) $\rho$ and  a cyclic resolution of identity $(P_0,\ldots,P_{d-1})$, every projection $P_k$ has the same weight, since for any fixed $k$ we have
$\rho(P_k)=\rho(\Phi^{*n}(P_k))$ for all $n$,
and, in particular for $n=k$, we get $\rho(P_k)=\rho(P_0)$.
\end{remark}

Now, we consider once again the special case of an OQRW $\mathfrak M$; with the notations introduced in previous sections, the associated unravelling is given by the operators $A_{i,j}=L_{i,j}\otimes\ketbra i j$, for $i,j\in V$.
\begin{prop}\label{prop_eqperiodicity}
A resolution of the identity $(P_0,\ldots,P_{d-1})$ is cyclic for an irreducible open quantum random walk $\M$ if and only if $P_k=\sum_{j\in V} P_{k,j}\otimes |j\rangle \langle j|$ for every $k$, with projectors $P_{k,j}$ satisfying the relation
\begin{equation}\label{eq_periodicity}
P_{k,i} L_{i,j} = L_{i,j} P_{k\submod 1,j}.
\end{equation}
\end{prop}

\pre

Assume that there exists an $A$-cyclic resolution of identity $(P_0,\ldots, P_{d-1})$.
Since $\M^*(P_k)=P_{k\submod 1}$, every $P_k$ is block-diagonal, \textit{i.e.} $P_k=\sum_j P_{k,j}\otimes |j\rangle \langle j|$, and from Proposition \ref{prop_Vcyclic}, for any $i,j$ in $V$: $ P_k \, L_{i,j}\otimes |i\rangle\langle j| =  L_{i,j}\otimes |i\rangle\langle j| \, P_{k\submod 1}$. This gives relation \eqref{eq_periodicity}. The converse is obvious.
\fin

\begin{remark}\label{remark-nonuniquedec}
For classical, irreducible, $d$-periodic Markov chains with stochastic matrix $K$, the cyclic components are uniquely determined and coincide with the irreducible communication classes $C_0,\ldots, C_{d-1}$ of the (aperiodic) Markov chain with transition matrix $K^d$. 
In the quantum context, the role of the partition $C_0,\ldots, C_{d-1}$, or, better yet, of the corresponding indicator functions $\ind_{C_0},\ldots, \ind_{C_{d-1}}$, is played by the cyclic projections $P_0,\ldots,P_{d-1}$. Indeed, notice that in the classical case $K \ind_{C_k}= \ind_{C_{k-1}}$ and, for the minimal dilation of this Markov chain, the cyclic projections $P_0,\ldots,P_{d-1}$ are uniquely determined as 
$P_k=\sum_{j\in C_k}|j\rangle\langle j|$.
However, an important difference should be underlined, with respect to the classical case:
in general, the resolution of the identity which verifies the definition of the period is not uniquely determined, since the decomposition of $\Phi^d$ into minimal irreducible components is not unique in general, as we will see in section \ref{section_stationary}.  An example of this fact can be easily constructed, as we now describe.
\end{remark}
\begin{example}
Take an OQRW $\mathfrak M$ with two sites and $\mathfrak h_1=\mathfrak h_2={\mathbb C}^2$, and introduce the matrix $R=\begin{pmatrix}0& -1\\ 1 & 0\end{pmatrix}$. Then we consider
$$
L_{11}=L_{21}=\frac{\i}{ \sqrt{2\,}}\, R,
\qquad
L_{12}=L_{22}= - \frac{\i}{ \sqrt{2\,}}\, R. 
$$
This $\mathfrak M$ is an irreducible OQRW (by a direct application of Proposition \ref{prop_ergodicityOQRW}) with period 2, and the cyclic projections $P_0,P_1$ can be chosen in different ways:
$$
P_0^{(x)} = |x\rangle\langle x| \otimes |1\rangle\langle 1| + |x\rangle\langle x| \otimes |2\rangle\langle2|, 
\quad
P_1^{(x)} = |Rx\rangle\langle Rx| \otimes |1\rangle\langle1| + |Rx\rangle\langle Rx| \otimes |2\rangle\langle2| 
$$
is a cyclic decomposition of the OQRW for any norm-one vector $x$ in $\cc^2$. As mentioned above, this is due to the fact that the map $\mathfrak M^2$ does not have a unique decomposition in irreducible components: $\M^2$ is the OQRW with all transition operators $L$ equal to $ \id /\sqrt 2$.  
\end{example}
We now discuss some results which will give us simple sufficient criteria for aperiodicity of an open quantum random walk.
\begin{lemme}\label{prop_caracaperiodicite}
Let $\M$ be a $d$-periodic open quantum random walk. Let $i,j\in V$ and $x\in\mathrm{Ran\,}P_{k,i}$,
 $y\in\mathrm{Ran\,}P_{k',j}$ for some $k,k'\in\{0,\ldots,d-1\}$. 
For any path~$\pi\in\mathcal P(i,j)$ of length $\ell$ one has
$ \langle y, L_\pi x\rangle =0$
unless $k'-k=\ell \ \mathrm{mod}\,d$.
\end{lemme}
\pre

Relation \eqref{eq_periodicity} implies that $L_\pi x$
belongs to the range of $P_{k\addmod\ell,j}$.
\fin

\begin{theo}\label{theo_caracaperiodicite}
Consider an irreducible open quantum random walk. For $i$ in~$V$, $x$ in $\mathfrak h_i$, define
\begin{equation}\label{def_Dix}
D(i,x)
= {\rm GCD}\{ \ell\ge 1, \exists \pi \in {\cal P}_\ell(i,i)\;{\rm s.t.}\; \langle x,L_\pi x\rangle \neq 0 \}.
\end{equation}
Then, for every $x$ in the range of $P_{k,i}$, the period $d$ is a divisor of $D(i,x)$.
In particular, if there exists $i$ in $V$ such that,
for all $x\in \mathfrak h_i,\, D(i,x)=1$,
then the open quantum random walk is aperiodic.
\end{theo}
\pre

Irreducibility implies that the defining set of $\ell$'s is nonempty, so that $D(i,x)$ is well-defined. The result follows from Lemma \ref{prop_caracaperiodicite}.
\fin

\begin{coro}\label{coro_perturbation}
Consider an irreducible open quantum random walk $\M$. If there exists $i$ in $V$ such that
\begin{equation}\label{eq_coroperturbation}
 \forall x\in \mathfrak h_i,\, \langle x,L_{i,i} x\rangle\neq0
\end{equation}
then the open quantum random walk is aperiodic.
\end{coro}

\begin{remark} 
The definition of the quantity $D(i,x)$ in Theorem~\ref{theo_caracaperiodicite} has an interpretation in terms of paths, and is reminiscent of the definition of the period for a state $i$ of a classical Markov chain with transition matrix $K$, \textit{i.e.} $D(i)=\mathrm{GCD}\,\{ \ell\geq 1\, |\, K_{ii}^{\ell}>0\}$. In addition, $D(i)$ coincides with \eqref{def_Dix} when applied to an OQRW which is a minimal dilation of the Markov chain. In the quantum context, however, $D(i,x)$ does not always coincide with the period, and, in particular, is not invariant with the argument $(i,x)$ even if the OQRW is irreducible (see Example \ref{example_Dnotconstant}). Even worse, the relation $d\, |\, D(i,x)$ may not hold if $x$ does not belong to the range of some $P_{i,k}$. Since the $P_k$ are \emph{a priori} unknown, the practical study of the period of an OQRW is difficult when simple sufficient conditions (such as the condition for aperiodicity given in Theorem \ref{theo_caracaperiodicite}) do not hold.
\end{remark}

\smallskip

In difficult cases, the following result can be helpful:
\begin{prop}\label{prop_caracaperiodic}
Consider an irreducible, finite, $d$-periodic open quantum random walk $\M$. If for some $i$ in $V$, and some $\ell$ prime with $d$, there exists a loop $\pi\in\mathcal P_\ell(i,i)$ of length $\ell$, such that $L_\pi$ is invertible, then $d$ is a divisor of $\mathrm{dim}\,\mathfrak h_i$.
\end{prop}

\pre

By Bezout's lemma, for any $k$ in $0,\ldots, d-1$ there exists an integer $a$ such that $a\ell=k \, \mathrm{mod}\, d$. Then $L_\pi^a P_{0,i} L_\pi^{-a}=P_{k,i}$, so that $\mathrm{dim}\, P_{k,i}$ does not depend on $k$. Therefore $\mathrm{dim}\, \mathfrak h_i= d\, \mathrm{dim}\, P_{0,i}$ and the conclusion follows. 
\fin

\begin{remark}\label{remark_loops} As a consequence of Corollary \ref{coro_perturbation}, starting from a finite irreducible periodic open quantum random walk $\mathfrak M$ we can perturb it into an aperiodic one, $\mathfrak M_{(\eps)}$, in different ways. If there exists $i_0$ in $V$ such that $L_{i_0,i_0}=0$ then one possible way is to define
\[ L_{i,j}^{(\eps)}=L_{i,j} \mbox{ if }j\neq i_0\quad \mbox{and} \quad L_{i,i_0}^{(\eps)}
=\left\{\begin{array}{ll}\sqrt{\eps\,}\,\mathrm{Id}& \mbox{ if }i=i_0,\\ \sqrt{1-\eps\,}\,L_{i,i_0} & \mbox{ if }i\neq i_0. \end{array}\right. \]  This is the analogue of ``adding a loop" for classical Markov chains. Another way is to ``add a loop" at every site, a method we will use in Example \ref{ex_Vn}. 
\end{remark}

For clarity we restate Proposition \ref{prop_cvgCPTPaperiodic} specifically for OQRWs:
\begin{theo}\label{theo_convergenceaperiodic}
Consider an irreducible, aperiodic and finite open quantum random walk $\M$.
For any state $\rho$, the sequence $(\mathfrak M^n(\rho))_n$
converges to the invariant state~$\rhoinv$ (which is unique, and faithful).
\end{theo}

\section{Ergodic properties of irreducible OQRWs}\label{section_recurrence}

We will now discuss the consequences of the previous theoretical results in terms of ergodic properties of irreducible open quantum random walks.
A first result in this direction is the following, which is a consequence of the ergodic theorem due to K\"ummerer and Maassen (\cite{KuMa}).
For completeness we give a self-contained proof in the present framework.
\begin{theo}[K\"ummerer-Maassen]\label{theo_kuma}
If the open quantum random walk $\mathfrak M$ is finite  then there exists a random variable $\rhoinv= \sum_{i\in V}\rhoinv(i)\otimes|i\rangle\langle i |$ with values in the set of invariant states on $\mathcal H=\bigoplus_{i\in V}\mathfrak h_i$ such that almost-surely,
\[ \frac 1n\sum_{k=0}^{n}\rho_k\otimes|X_k\rangle\langle X_k| \underset{n\to\infty}{\longrightarrow} \sum_{i\in V}\,\rhoinv(i)\otimes|i\rangle\langle i |.\]
\end{theo}
\pre

Let $\eta_n$ be the state $\rho_n \otimes |X_n\rangle\langle X_n|$. Denote by $\mathcal F_n$ the $\sigma$-algebra generated by $\eta_k$ for $k\leq n$, and let
\[m_n = \sum_{k=0}^n \eta_k- \sum_{k=0}^{n-1} \mathfrak{M}( \eta_{k}).\]
We have, from \eqref{eq_exprhon} above, $\mathbb E(m_{n+1} - m_n | \mathcal F_n) =0$ so that $(m_n)_{n}$ is a martingale, and since $\|m_{n+1}-m_n\|=\| \eta_{n+1} - \mathfrak M (\eta_n)\|$
is uniformly bounded, we can apply the law of large numbers for  martingales with uniformly bounded increments. Therefore,  
$
 \frac{1}{n} \sum_{k=0}^{n}\eta_k - \frac{1}{n} \sum_{i=0}^{n-1} \mathfrak{M} (\eta_k) \to 0 
$
where convergence is meant in the almost-sure sense. In turn, this implies for any $N\in \nn^*$,
 \[
 \frac{1}{n} \sum_{k=0}^{n}\eta_k - \frac{1}{n} \sum_{k=0}^{n-1} \mathfrak{M}^N (\eta_k) \to 0 
 \]
 so that
  \[
 \frac{1}{n} \sum_{k=0}^{n}\eta_k - \frac{1}{n} \sum_{k=0}^{n-1} \frac{\id+\mathfrak{M}+\ldots+\mathfrak{M}^{N-1}}N (\eta_k) \to 0.
 \]
For any state $\eta$, $ \frac{\id+\mathfrak{M}+\ldots+\mathfrak{M}^{N-1}}N (\eta)$ converges when $N$ goes to infinity to an invariant state. This can be seen viewing $\M$ as a contraction on the Hilbert-Schmidt space $\mathcal I_2(\mathcal H)$, \textit{i.e.} $\mathcal B(\mathcal H)$ equipped with the scalar product $\tr(A^*B)$. This invariant state must be of the form $P\eta$, where $P$ is a linear operator on $\mathcal I_1(\mathcal H)$. The operator $P$ can be approximated uniformly by $\frac{I + \mathfrak{M} + \ldots + \mathfrak{M}^{N-1}}{N}$, therefore 
$
 \frac{1}{n} \sum_{k=0}^{n} \eta_k - \frac{1}{n} \sum_{k=0}^{n-1} P \eta_k\to 0.
$ On the other hand $P\mathfrak M = P$ implies that $ \ee(P\eta_{n+1}|\mathcal F_n)=P\eta_n$, \textit{i.e.} $(P\eta_n)_n$ is a bounded martingale, so $\frac{1}{n} \sum_{k=0}^{n} P \eta_k$ converges almost-surely to some invariant state. This concludes the proof.
\fin

A direct consequence of Theorem \ref{theo_kuma} (of which we shall preserve the notations) and of our previous observations on the form of $\rhoinv$ is the following:
\begin{coro}\label{coro_rec1}
If the open quantum random walk $\M$ is finite and irreducible with invariant (and faithful) state
$\rhoinv=\sum_{i\in V}\rhoinv(i)\otimes|i\rangle\langle i |,$
then for all  $i$~in~$V$, define 
$ N_n(i)=\mathrm{card}\,\{k\leq n\, |\, X_k=i\}$. We have
\begin{eqnarray*}
 \frac{N_n(i)}{n} &\underset{n\to\infty}{\longrightarrow}& \tr\,\rhoinv(i)\mbox{ almost-surely,}\\
\frac1n\sum_{k=0}^{n-1}\pp(X_k=i) &\underset{n\to\infty}{\longrightarrow}& \tr\,\rhoinv(i),\\
\frac 1{N_n(i)}\sum_{k=0}^{n-1}\rho_k \, \ind_{X_k=i}&\underset{n\to\infty}{\longrightarrow}& \frac{\rhoinv(i)}{\tr\,\rhoinv(i)}\mbox{ almost-surely.}
\end{eqnarray*}
\end{coro}
\pre

This is simply obtained by examination from Theorem \ref{theo_kuma}.
\fin

\begin{remark} The only new result that we brought to the above picture is Theorem~\ref{theo_unicity}, which tells us that the state $\rhoinv$ is unique and faithful, and in particular $\tr\,\rhoinv(i)>0$ for any $i$ in $V$. This implies that, for any irreducible open quantum random walk with an invariant state $\rhoinv$, one has
$$ \mbox{for all }i\in V, \quad \pp(X_n=i\ \mbox{infinitely often})=1,$$
$$ \mbox{for all }i\in V,\, x\in \mathfrak h_i,\quad \pp(\langle x,\rho_n(i)\,x\rangle\,\ind_{X_n=i}>0 \ \mbox{infinitely often})=1.$$ 
The first statement has an immediate interpretation in terms of ``spatial recurrence" (every site $i$ in $V$ is visited infinitely often), the second one is stronger and can be seen as ``spatial and internal recurrence".
\end{remark}
\smallskip

The second ergodic result of this section is a consequence of Theorem \ref{theo_ergodicconvergence}.
\begin{coro}\label{coro_rec2}
If the open quantum random walk $\M$ is irreducible with invariant (and faithful) state
$\rho_{\mathrm{inv}}=\sum_{i\in V}\rhoinv(i)\otimes|i\rangle\langle i |,$
then for all $i$ in~$V$,
\begin{eqnarray*}
\frac1n \sum_{k=0}^{n-1} \pp(Q_k=i)&\underset{n\to\infty}{\longrightarrow}& \tr\rhoinv(i)\\
\frac 1n\sum_{k=0}^{n-1}\mathfrak M ^k(\rho,i) &\underset{n\to\infty}{\longrightarrow}& \rhoinv(i) \mbox{ in the weak-$\ast$ sense}.
\end{eqnarray*}
\end{coro}

\begin{remark} The assumption that there exists an invariant state is necessary in Corollary \ref{coro_rec2} (contrary to Corollaries \ref{coro_rec1} and \ref{coro_rec3}, where it is always true and only stated to establish notations), because we do not assume finiteness of $\M$. Since for all $n\in \nn$ and $i$ in $V$ one has $\pp(X_n=i)=\pp(Q_n=i)$, the first statement of Corollary \ref{coro_rec2} is a refinement of the second statement of Corollary~\ref{coro_rec1}.
\end{remark}

Our third corollary is a consequence of Theorem \ref{theo_convergenceaperiodic}. It is an improvement of the previous result in the case where we have aperiodicity.

\begin{coro}\label{coro_rec3}
If the open quantum random walk $\M$ is finite, irreducible and aperiodic with invariant (and faithful) state
$\rho_{\mathrm{inv}}=\sum_{i\in V}\rhoinv(i)\otimes|i\rangle\langle i |,$
then for all $i$ in $V$,
\begin{eqnarray*}\pp(Q_n=i) &\underset{n\to\infty}{\longrightarrow}& \tr\rhoinv(i),\\
 \mathfrak M ^n(\rho,i) &\underset{n\to\infty}{\longrightarrow}& \rhoinv(i).\end{eqnarray*}
\end{coro}

\begin{remark}\ 
\vspace{-0.5em}\begin{itemize}
\item Corollary \ref{coro_rec2} seems rather useless, from an operational point of view: there is no joint realization of the different $Q_k$ or $\mathfrak M^k(\rho,i)$ for different $k$; or, in other terms, measuring $Q_k$ disrupts the existence of $Q_{k'}$ or $\mathfrak M^{k'}(\rho,i)$ for ${k'}>k$. Corollary \ref{coro_rec3}, on the other hand, is operational, and tells us that, if the system is left to evolve for a large time, then a single measurement will give the position $i$ with approximate probability $\tr\rhoinv(i)$, and the state after that unique measurement will be approximately $\frac{\rhoinv(i)}{\tr\rhoinv(i)}$. These limiting quantities are the same as those that appear for limits \emph{with} measurements. This results display evident similarities with the behavior of classical Markov chains.
\item The results concerning the probabilities may not seem surprising, precisely because of the equality of the laws of $X_n$ and $Q_n$ and of the K\"ummerer-Maassen theorem; the convergence in Corollary \ref{coro_rec3}, however, is completely new. Results regarding the induced states $\M^n(\rho,i)$ in Corollary \ref{coro_rec2} are new and show that the limits are the same whether we do an infinite number of measurements (for the sequence $(\rho_n)_n$), or just one but after an ``infinite" time (for $(\M^n(\rho))_n$). In addition, all convergences in Corollary~\ref{coro_rec3} (without averaging) are new.
\item Example \ref{ex_Vn} shows that irreducibility is a necessary assumption for Corollary \ref{coro_rec3}, as one could expect from analogous results for classical Markov chains.
\end{itemize}
\end{remark}

\section{Reducible OQRWs and communication classes} \label{section_nonirreducible}

In this section we study the failure of irreducibility for an open quantum random walk. Considering reducible OQRWs, the first natural problem one has to face is how to characterize reducibility and how to determine reducing, possibly minimal, components.

A reasonable way to proceed, mimicking what happens for classical Markov chains, is to define a communication relation between vectors of the Hilbert space~${\cal H}$: this relation should be an equivalence relation constructed in such a way that the induced equivalence classes are the irreducible components of the map $\M$. We will see that it is possible to do this in a way which is consistent with the classical case.

However, it is important to immediately underline that the quantum case displays peculiar features: the decomposition of the Hilbert space ${\cal H}$ as the direct sum of irreducible components is not unique in general. This is not at all surprising if one thinks about the structure of invariant states for a CP-TP map (see \cite{BN}, from which we take much of our inspiration): essentially, the quantum peculiarity is related to the fact that there can exist stationary states which are not simple convex combinations of the stationary states on each irreducible component.

\smallskip

 Following Baumgartner and Narnhofer (\cite{BN}), we call a closed vector space~$\mathcal V$ an \textit{enclosure} (for $\M$) if $\mathrm{supp} \, \rho\subset \mathcal V$ for $\rho$ a positive, trace-class operator implies $\mathrm{supp} \, \mathfrak M(\rho)\subset \mathcal V$. The next proposition will be extremely useful. The fact that the support of an invariant state is an enclosure is however a known result, see \textit{e.g.} section 1 in~\cite{FV}. 

From now on we fix an OQRW $\M$ with the same notation as in section \ref{section_OQRWs}.
 
\begin{prop}\label{lemma_enctraj}
\ \\
\vspace{-1.8em}
\begin{enumerate}
\item If $\mathcal V$ is a closed subspace of $\mathcal H$, then it is an enclosure if and only if it is stable by $L_\pi\otimes |j\rangle\langle i |$ for any $i,j$ in $V$ and $\pi\in\mathcal P(i,j)$. In particular, if a vector $x=\sum_{i\in V}{x_i}\otimes |i\rangle$ is in $\mathcal V$ then $ {L_\pi x_i}\otimes\vec{j}\in\mathcal V$.
\item A projection $P$ reduces $\mathfrak M$ if and only if 
$$
P(L_\pi\otimes |j\rangle\langle i |)P=(L_\pi\otimes |j\rangle\langle i |)P
\quad \mbox{for all $i,j$ in $V$ and } \pi\in\mathcal P(i,j).
$$
\item
If $\mathcal V$ is an enclosure, then $(L_\pi\otimes |j\rangle\langle i |)(\mathcal V)\subset \mathfrak h_j\cap \mathcal V$  for any $i,j$ in $V$ and~$\pi\in\mathcal P(i,j)$, and $\bigoplus_{j\in V} {\rm Vect}\{(L_\pi\otimes |j\rangle\langle i |)(\mathcal V),i\in V,\pi \in\mathcal P(i,j)\}$ is also an enclosure.

\item The support of an invariant state is an enclosure.
\end{enumerate}
\end{prop}

\pre

\begin{enumerate}
\item Suppose that $\mathcal V$ is an enclosure. Remark that, for any positive integer $\ell$ and any $x=\sum_{i\in V}{x_i}\otimes |i\rangle$ in $\mathcal H$, one has
\begin{equation}
\label{eq_Mell} \M^\ell (\ketbra xx) = \sum_{i,j\in \, V} \sum_{\pi\in\mathcal P_\ell(i,j)} \ketbra{L_\pi x_i}{L_\pi x_i}\otimes \ketbra jj;
\end{equation}
so, if $x$ is in $\mathcal V$, then every ${L_\pi x_i}\otimes\vec{j}$ is in $\mathcal V$. 
Conversely, let us now suppose that $\mathcal V$ is stable under the action of the operators $L_\pi\otimes |j\rangle\langle i |$. Starting from~$\rho$ with $\mathrm{supp}\, \rho$ in $\mathcal V$, then considering its spectral decomposition and~\eqref{eq_Mell} above shows that $\mathrm{supp}\,\M^\ell(\rho)\subset \mathcal V$. This proves the first statement.
\item Just recall that a subspace of $\mathcal H$ is the support of a reducing projection if and only if it is an enclosure. This point is then an immediate consequence of the previous one. 
\item This point also follows immediately from $1$.

\item Consider an invariant state $\rho_0$ and a state $\rho$ with support contained in $\mathrm{supp}\,\rho_0$. Then there exists a weak approximation of $\rho$ by an increasing sequence of finite-dimensional operators $\rho_n$ with $\mathrm{supp}\,\rho_n\subset\mathrm{supp}\,\rho_0$. Furthermore, for every $n$, there exists a $\lambda_n$ such that $\rho_n\leq \lambda_n\rho_0$, so that $\M(\rho_n)\leq \lambda_n \rho_0$ and $\mathrm{supp}\,{\mathfrak M}(\rho_n)\subset\mathrm{supp}\,\rho_0$. The sequence $\M(\rho_n)$ is increasing and weakly convergent to $\M(\rho)$ so that $\mathrm{supp}\,\M(\rho)\subset\mathrm{supp}\, \rho_0$, which proves that $\mathrm{supp}\, \rho_0$ is an enclosure.\ \fin
\end{enumerate}
 
In general, it is not true that all the reducing projections are diagonal, \textit{i.e.} of the form $\sum_{i\in V}P_i\otimes\ketbra ii$, but by the previous Proposition, point $3$, if $\mathfrak M$ is reducible, then it admits at least one block-diagonal reducing projection. So the reducibility of an OQRW can be established considering only block-diagonal projections.
Moreover, notice that the support projection of an invariant state is block-diagonal, \textit{i.e.} of the form
$P=\sum_{i\in V} P_i\otimes \ketbra ii$.

We can characterize the block-diagonal projections reducing $\mathfrak M$ using the unravelling of $\mathfrak M$.

\begin{prop}\label{SubhProj}
An orthogonal block-diagonal projection $P=\sum_j P_j\otimes|j\rangle\langle j|$ reduces $\mathfrak M$ if and only if $\mathrm{Ran}\,L_{i,j}P_j\subset \mathrm{Ran}\,P_i$,
(i.e. $L_{i,j}P_j=P_i L_{i,j}P_j$) for all $i$ and $j$ in $V$.
Equivalently, a closed subspace of the form $\mathcal V=\oplus_{i\in V}\mathcal V_i$, with~$\mathcal V_i\subset \mathfrak h_i$, is an enclosure if and only if $L_{i,j} \mathcal V_j \subset \mathcal V_i$ for all $i$ and $j$.
\end{prop}
\pre 

It is clear that it is sufficient to prove only the first statement, due to the relation between enclosures and reducing projections that we have recalled also in the previous proof.
So, take a reducing projection $P=\sum_{j\in V} P_j\otimes|j\rangle\langle j|$.
By point 2 in the previous proposition, it is necessary that the range of $P$ is invariant under the action of all operators of the form
$L_{i,j}\otimes |i\rangle\langle j |$ and so the relation $L_{i,j}P_j=P_i L_{i,j}P_j$ for all $i$ and $j$ in $V$ immediately follows. 
The reverse implication is also easy to obtain using again the characterization in point 2 of previous proposition and the fact that any operator $L_\pi\otimes |i\rangle\langle j |$, for a path~$\pi=(i_0,i_1,....i_{\ell})\in {\mathcal P}(i,j)$ ($i=i_0,\, j=i_\ell$) of length $\ell$, is the composition of the operators $L_{i_{k+1},i_k}\otimes |i_{k+1}\rangle\langle i_k |$, with $k=0,...\ell-1$.
 \fin

\begin{coro}
When $P=\sum_{j\in V} P_j\otimes|j\rangle\langle j|$ is a reducing projection, each $P_j$ is a projection on a subspace preserved by $L_{jj}$.
\end{coro}

\begin{remark}\label{InvariantStates}
Suppose that, for all sites $i$ and $j$, there exists ``a path of invertible operators which connects them'', i.e. a path $\pi\in{\cal P}(i,j)$ such that $L_\pi$ is invertible. In this case, the previous proposition proves that, if $P=\sum_{j\in V} P_j\otimes|j\rangle\langle j|$ is a reducing projection, then $\mathrm{rank}\,P_i=\mathrm{rank}\, P_j$ for any $i,j$ in $V$. In particular, if a state $\rho=\sum_{i\in V}\rho_i\otimes \ketbra ii$ is invariant, then $\rho_i\neq 0$ for all $if$ (i.e. an invariant state is supported by all sites), and if $\rho_i$ is faithful on $\mathfrak h_i$ for some index $i$, then~$\rho_j$ is faithful on $\mathfrak h_j$ for any $j$ in $V$.
\end{remark}

\smallskip
The next notion, of enclosure generated by a single vector in $\mathcal H$, will be crucial in our analysis of decompositions of reducible OQRWs:
\begin{defi}
For $\phi$ in $\mathcal H$, we denote by $\enc(\phi)$ the closed vector space 
$$\enc(\phi)=\overline{ \mathrm{Vect}\,\bigcup_{i,j\in V}\{(L_\pi\otimes \ketbra ji) \, \phi \, |\, \pi\in\mathcal P(i,j)\}}.$$
Consistently with Proposition \ref{lemma_enctraj}, we will consider specifically enclosures of vectors $x\otimes \vec i$, which take the form
$$\enc(x\otimes \vec i)=\overline{ \mathrm{Vect}\,\bigcup_{j\in V}\{{L_\pi\, x}\otimes|j\rangle   \, |\, \pi\in\mathcal P(i,j)\}}.$$
\end{defi}
We will be mostly interested in enclosures that are minimal but non-trivial (\emph{i.e.} not equal to $\{0\}$). From now on, the term \emph{minimal enclosure} will refer to minimal, non-trivial enclosures. 
The following lemma contains relevant properties of enclosures. 
\begin{lemme}\label{lemma_lemenc}\ \vspace{-0.5em}
\begin{itemize}
\item The space $\enc( x\otimes\vec i)$ is the smallest enclosure containing $x\otimes \vec i$.
\item Any minimal enclosure is of the form $\enc(x\otimes\vec i)$.
\item If two minimal enclosures $\enc(x\otimes\vec i)$ and $\enc (y\otimes\vec j)$ are distinct then they are in direct sum.
\end{itemize}
\end{lemme}

\pre

All statements follow from Proposition \ref{lemma_enctraj}.
\fin

\begin{remark}\label{remark_encxi}
In the same way that the specific form of $\M(\rho)$ led us to consider only states $\rho$ of the form $\rho=\sum_{i\in V}\rho(i) \otimes |i\rangle \langle i|$, Proposition \ref{lemma_enctraj} shows that vectors of the form $x\otimes \vec i$ are of particular interest. In particular, any minimal enclosure will be generated by a vector $x\otimes \vec i$.
It is not true, however, that any $\enc (x\otimes \vec i)$ is a minimal enclosure, as the following example shows.
\end{remark} 
 \begin{example}\label{ex_2}
Take $V=\{1,2,3\}$ with
$$ L_{1,2}=L_{2,3}=L_{3,1}=\frac1{\sqrt 5}\begin{pmatrix}1 & 0\\ 0 & 2\end{pmatrix}  \qquad L_{2,1}=L_{3,2}=L_{1,3}=\frac1{\sqrt 5}\begin{pmatrix}2 & 0\\ 0 & 1\end{pmatrix}.$$
One can see that for $k=1,2$, $\enc(e_k\otimes |1\rangle)$, are minimal enclosures, equal to~$\cc\, e_k\otimes\ell^2(V)$, but the space $\enc((e_1+e_2)\otimes |1\rangle)$ is equal to $\cc^2\otimes\ell^2(V)$. \end{example}

\begin{remark} \label{remark-siteconnection}
Let us return to the notion of irreducibility, as introduced in Definition~\ref{defi_irroqrw}: an open quantum random walk $\M$ is irreducible if for any $i,j$ in $V$, one has $i\connect j$, which by Proposition \ref{prop_ergodicityOQRW} is defined by the equivalent conditions
\begin{equation}\label{eq_connect1} \forall \, x\in\mathfrak h_i, \, y\in \mathfrak h_j,\ \exists \pi\in\mathcal P(i,j) \mbox{ such that }\langle{y},{L_\pi x}\rangle\neq 0,
\end{equation}
\begin{equation}\label{eq_connect2}
\forall \, x\in\mathfrak h_i, \, y\in \mathfrak h_j,\ y\in \overline{\mathrm{Vect}\,\{L_\pi\, x\, |\, \pi\in \mathcal P(i,j)\}}.
\end{equation}
From the above discussion it is clear that both conditions can be characterized using enclosures:
\begin{eqnarray*}
\exists \pi\in\mathcal P(i,j) \mbox{ such that }\langle{y},{L_\pi x}\rangle\neq 0 &\Leftrightarrow&  y \not\in \enc(x\otimes \vec i)^\perp\\
y\in \overline{\mathrm{Vect}\,\{L_\pi\, x\, |\, \pi\in \mathcal P(i,j)\}} &\Leftrightarrow&  y \in \enc(x\otimes\vec i).
 \end{eqnarray*}
As we will see below, in Proposition \ref{prop_coherence}, the orthogonal of an enclosure can be related to another enclosure. This will allow us to strengthen the connection between the two notions above.
\end{remark}

The above discussion gives immediately:
\begin{lemme}
An open quantum random walk $\M$ is irreducible if and only if~$\mathcal H$ is a minimal enclosure, or equivalently, if $\mathcal H=\enc(x\otimes \vec i)$ for any $x\otimes \vec i$ in~$\mathcal H$.
\end{lemme}

To emphasize the picturesque aspect of our definition of irreducibility, we define the following notion of accessibility among vectors, denoted by $\connect$. We remark that the notation, $\connect$, is the same we used in Definition \ref{defi_irroqrw}, but this should not generate confusion, the difference being clear in the arguments we use: in the previous case, the connection $\connect$ is between sites $i$, $j$ in $V$, whereas here it is between vectors $\phi$, $\psi$ of the Hilbert space $\mathcal H$.
\begin{defi}
For $\phi$, $\psi$ in $\mathcal H$, we denote $\phi \connect \psi$ if $\psi\in\enc(\phi)$, and $\phi \biconnect \psi$ if $\phi \connect \psi$ and $\psi \connect \phi$.
\end{defi}
Again, we will be specifically interested in the relation $\connect$ between vectors of the form $x\otimes \vec i$, $y\otimes \vec j$ and we have immediately
\[
x\otimes \vec i \connect y\otimes \vec j \Leftrightarrow y \in \overline{\{L_\pi \, x \, |\, \pi \in \mathcal P(i,j)\}}.
\]

Going back to the connection between sites, we can also add that, due to Remark \ref{remark-siteconnection},
$$
i\connect j 
\qquad\Leftrightarrow \qquad
x\otimes \vec i \connect y\otimes \vec j \mbox{ for all } x\in\mathfrak h_i, y\in\mathfrak h_j.
$$

The following Proposition can easily be proven:
\begin{prop}
The relation $\connect$ on $\mathcal H$ is transitive, and $\biconnect$ is an equivalence relation. Any minimal enclosure is an equivalence class of $\biconnect$.
\end{prop} 
\begin{remark}
Every equivalence class of a vector $x\otimes \vec i $ by $\biconnect$ is a subset of~$\mathcal H$ contained in $\enc (x\otimes \vec i)$, but it may fail to be an enclosure and even a subspace. A minimal dilation of a classical Markov chain with a proper transient class easily gives an example of an equivalence class that is not an enclosure. For an example where an equivalence class is not a subspace, consider Example \ref{ex_notsubspace} below.
\end{remark}
\begin{example}\label{ex_notsubspace}
Consider $V=\{1,2\}$, $\h_1=\h_2=\cc^2$ with canonical basis denoted by $(e_1,e_2)$, and introduce the OQRW $\M$ with transitions
\[L_{1,1}=L_{2,2}=\frac1{\sqrt 2}\begin{pmatrix}0&1\\1&0\end{pmatrix} \qquad L_{1,1}=L_{2,2}=\frac1{\sqrt 2}\,\id.\]
Then the only minimal enclosures are
\[E_+=\enc\big((e_1+e_2)\otimes \vec 1 \big) = \cc\, (e_1+e_2)\otimes\ell^2(V)\]
\[E_-=\enc\big((e_1-e_2)\otimes \vec 1 \big) = \cc\, (e_1-e_2)\otimes\ell^2(V)\]\end{example}
and for any $x\not\in \cc(e_1+e_2)\cup \cc(e_1-e_2)$ one has
\[\enc(x\otimes \vec 1) = \enc(x\otimes \vec 2) =  \H.\]
Therefore, for such an $x$, the equivalence class of $x\otimes \vec 1$ is $\H \,\setminus (E_+\cup E_-)$. 
\section{Decompositions of OQRWs and invariant states}\label{section_stationary}

In this section we wish to focus on the behavior of an OQRW $\mathfrak M$ on the so-called fast recurrent subspace, i.e. the support of the $\mathfrak M$-invariant states.
We decompose the corresponding restriction of $\mathfrak M$ into a ``direct sum" of irreducible OQRWs $\mathfrak M_k$, 
establish when this decomposition is unique, and study how the different irreducible components interact.
We follow the lines traced in \cite{BN} for quantum evolutions on finite dimensional Hilbert spaces;
we will state and prove generalizations to infinite dimension. As we will see in Proposition \ref{prop_coherence}, the form of invariant states is dictated by the unicity or non-unicity of the decompositions into minimal enclosures, and Lemma \ref{lemma_unicitydec} shows that non-unicity is related to the existence of mutually non-orthogonal minimal enclosures. 

To further study the stationary states of an OQRW, we recall some notation. Inspired by \cite{FV}, we denote: 
\[ \mathcal R = \mathrm{sup}\{\mathrm{supp}\,\rho\, |\, \rho \mbox{ an invariant state}\} . \]
This space is often called the fast recurrent space.

\begin{remark}\label{remark_RDfinite}
The above definition of $\mathcal R$ is unfortunately not explicit, and makes a (small) part of Theorem \ref{theo_invariantstates} describing stationary states tautological. In the finite dimensional case, $\mathcal R$ can be equivalently described (as is done in \cite{BN}) without reference to the set of invariant states, as $\mathcal R=\mathcal D^\perp$, where $\mathcal D$ is defined~by
\[\mathcal D = \{\phi\in \mathcal H\, |\,\langle\phi,\M^n(\rho)\,\phi\rangle\underset{n\to\infty}{\longrightarrow}0 \mbox{ for any state }\rho\}.\]
\end{remark}

The following Lemma is an immediate consequence of Proposition \ref{lemma_enctraj}.
\begin{lemme}
The subspace $\mathcal R$ is an enclosure.
\end{lemme}

We let $\mathcal D =\mathcal R^\perp$, which is characterized as
\[\mathcal D = \{\phi \in \mathcal H \, |\, \langle \phi,\rho\,\phi\rangle =0 \mbox{ for any invariant state } \rho\}.\] From the block-diagonal structure of $\mathfrak M^n (\rho)$ for $\rho$ any state, we clearly have
\[ \mathcal R = \bigoplus_{i\in V} \mathcal R_i \mbox{ with } \mathcal R _i \subset \mathfrak h_i,
 \qquad 
 \mathcal D = \bigoplus_{i\in V} \mathcal D_i \mbox{ with } \mathcal D _i \subset \mathfrak h_i. \]

Since our main interest is to investigate the invariant states of open quantum random walk $\M$ on $\mathcal H$, we will be interested in decomposing $\mathcal R$, not~$\mathcal D$, into irreducible subsystems.
\smallskip

We will use the following results, which were stated in \cite{BN} in the finite dimensional case. We extend them here to infinite dimension.
\begin{prop}\label{prop_coherence}
If $\mathcal V$ and $\mathcal W$ are two subspaces of $\mathcal H$ such that $\mathcal V \cap \mathcal W =\{0\}$ and $\rho$ is a state with support in $\mathcal V \oplus \mathcal W$, then denote
\[\rho_{\mathcal V}= P_{\mathcal V}\, \rho \,P_{\mathcal V},\quad \rho_{\mathcal W}= P_{\mathcal W}\, \rho \,P_{\mathcal W}\, \quad \rho_{\mathcal C}= P_{\mathcal V}\, \rho\, P_{\mathcal W}, \quad \rho_{\mathcal C}'=P_{\mathcal W}\,\rho \,P_{\mathcal V} \]
so that $\rho = \rho_{\mathcal V} + \rho_{\mathcal W} + \rho_{\mathcal C}+\rho_{\mathcal C}'$. Decompose $\M(\rho)$ in a similar way.
\begin{enumerate}
\item If $\mathcal V$ is an enclosure, then $P_{\mathcal W}\, \M (\rho_C+\rho_C')\, P_{\mathcal W}=0$.
\item If $\mathcal V$ is an enclosure, then so is $\mathcal V^\perp\cap\mathcal R$.
\item If $\mathcal V$ and $\mathcal W$ are enclosures, then
\[ 
\M (\rho)_{\mathcal V} =\M(\rho_{\mathcal V}) \quad \M (\rho)_{\mathcal W} =\M(\rho_{\mathcal W})
\quad \M (\rho)_{\mathcal C} =\M(\rho_{\mathcal C}) \quad \M (\rho)_{\mathcal C}' =\M(\rho_{\mathcal C}').\]
\item  A subspace of $\mathcal R$ is a minimal enclosure if and only if it is the support of an extremal invariant state. In particular, if $\mathcal V\subset\mathcal R$ is an enclosure, then it contains a (non-trivial) minimal enclosure.
\item If $\rho$ is $\M$-invariant and $\mathcal V$ and $\mathcal W$ are two minimal enclosures contained in~$\mathcal R$, such that the decomposition of $\mathcal V \oplus \mathcal W$ into a sum of minimal enclosures is unique, then $\rho_{\mathcal C}=0$ and $\rho_{\mathcal C}'=0$.
\end{enumerate}
\end{prop}
\pre

We essentially borrow the main ideas of the proofs from \cite{BN}, adding some variations when required by the infinite dimensional setting. 
\begin{enumerate}
\item To prove the first point we define
$\kappa_{\pm \eps}=\frac1\eps \,\rho_{\mathcal V}\pm \,(\rho_{\mathcal C}+\rho_{\mathcal C}')+\eps \, \rho_{\mathcal W}.$
We have~$\kappa_{\pm\eps}\geq 0$ (as can be checked from $\langle u, \kappa_{\pm\eps}\, u \rangle = \langle u_{\pm\eps}, \rho\, u_{\pm\eps} \rangle$ where $ u_{\pm\eps}=\frac1{\sqrt \eps}\,P_{\mathcal V}u + \sqrt \eps\, P_{\mathcal W}u$), so that $\mathfrak M(\kappa_{\pm\eps})\geq 0$, and, because $\mathcal V$ is an enclosure, the support of $\mathfrak M(\rho_{\mathcal V})$ is contained in $\mathcal V$, so that
\[ P_{\mathcal W}\,\mathfrak M(\kappa_{\pm \eps})\, P_{\mathcal W} = \pm P_{\mathcal W}\,\mathfrak M(\rho_{\mathcal C}+\rho_{\mathcal C}')\,P_{\mathcal W}+ \eps \,P_{\mathcal W}\, \mathfrak M (\rho_{\mathcal W}) \, P_{\mathcal W}.\]
This must be $\geq 0$ for any $\eps$, and by necessity $ P_{\mathcal W}\,\mathfrak M(\rho_{\mathcal C}+\rho_{\mathcal C}')\,P_{\mathcal W}=0$. 

\item Consider $\mathcal W=\mathcal V^\perp$ and $\eta$ any invariant state; then 
\[\eta_{\mathcal V}+\eta_{\mathcal W}+\eta_{\mathcal C}+\eta_{\mathcal C}' = \M(\eta_{\mathcal V}) + \M(\eta_{\mathcal W})+\M(\eta_{\mathcal C})+\M(\eta_{\mathcal C}').\]
Projecting by $P_{\mathcal W}$ this yields $\eta_{\mathcal W}= P_{\mathcal W} \M(\eta_{\mathcal W})P_{\mathcal W} $, so that $P_{\mathcal V}\,\M(\eta_{\mathcal W})\, P_{\mathcal V}$ is positive with zero trace. Therefore $P_{\mathcal V}\,\M(\eta_{\mathcal W})\, P_{\mathcal V}=0$ which implies $P_{\mathcal V}\,\M(\eta_{\mathcal W})=\M(\eta_{\mathcal W})\, P_{\mathcal V}=0$ and so $\eta_{\mathcal W}=\M(\eta_{\mathcal W})$. 
As the support of a stationary state, $\mathrm{supp}\,\eta_{\mathcal W}= \mathrm{supp}\,\eta\cap\mathcal V^\perp$ is an enclosure. Taking the supremum over all possible invariant states $\eta$, this tells us that $\mathcal R \cap \mathcal V^\perp$ is also an enclosure.

\item If both $\mathcal V$ and $\mathcal W$ are enclosures, then by point 1, and the fact that $\mathrm{supp}\,\M(\rho_{\mathcal V})\subset \mathcal V$ and $\mathrm{supp}\,\M(\rho_{\mathcal W})\subset \mathcal W$, we have 
\begin{equation}\label{eq_etapeinvariance}
\M(\rho_{\mathcal C})+\M(\rho_{\mathcal C}')=\M(\rho)_{\mathcal C} +\M(\rho)_{\mathcal C}'.
\end{equation}
Now remark that if \textit{e.g.} $\phi\in \mathcal V$ and $\psi\in \mathcal W$, then for any $i$ and $j$ in $V$ we have
\[\big(L_{i,j}\otimes\ketbra ij\big) \phi \in \mathcal V\quad\mbox{ and }\big(L_{i,j}\otimes\ketbra ij\big)\psi \in \mathcal W. \]
Therefore, \eqref{eq_etapeinvariance} actually implies $\M(\rho_{\mathcal C})=\M(\rho)_{\mathcal C}$ and $\M(\rho_{\mathcal C}')=\M(\rho)_{\mathcal C}'$.
\item If $\mathcal V$ is a minimal enclosure contained in $\mathcal R$, then there exists an $\M$-invariant state $\rho$ such that $\rho_{\mathcal V}=P_{\mathcal V}\rho P_{\mathcal V} \neq 0$. By point $3$, we have 
$\rho_{\mathcal V}=\M(\rho)_{\mathcal V}=\M(\rho_{\mathcal V})$, and so
 $\rho_{\mathcal V}$ is (up to normalization) an invariant state of $\M_{|\mathcal I_1(\mathcal V)}$. Since ${\mathcal V}$ is irreducible, by Theorem \ref{theo_unicity}, $\M_{|\mathcal I_1(\mathcal V)}$ has a unique invariant state, which has support equal to $\mathcal V$. Therefore, $\rho_{\mathcal V}$ is a state with support $\mathcal V$. This $\rho_{\mathcal V}$ must be extremal since $\rho_{\mathcal V}=t\, \rho_1 + (1-t)\, \rho_2$ with~$\rho_1$, $\rho_2$ invariant states and $t\in]0,1[$ would imply that $\rho_1$, $\rho_2$ are invariant states with support in $\mathcal V$ but then by unicity, $\rho_{\mathcal V}=\rho_1=\rho_2$.

Conversely, if $\mathcal V=  \mathrm{supp}\, \rho$ with $\rho$ an extremal invariant state, then $\mathcal V$ must be an enclosure.
If, by contradiction, we suppose it is not minimal, then there exists an enclosure $\mathcal W$ with $\mathcal W \subset\mathcal V\subset \mathcal R$; then, using point~$2$ and repeating the arguments of the previous implication would yield  the existence of two states (up to normalization) $\rho_{\mathcal W}$ and $\rho_{\mathcal W^\perp\cap \mathcal V}$ which are  invariant, of which $\rho$ is a convex combination. The extremality of $\rho$ implies that $\mathcal W$ is either $\{0\}$ or $\mathcal V$ and so $\mathcal V$ is minimal.

To prove the last statement, observe that by definition there exists an invariant $\rho$ such that $\mathcal V\cap\, \mathrm{supp}\,\rho\neq\{0\}.$ By point 3, $\mathcal V$ contains the support of the invariant state $\rho_{\mathcal V}$. By the Krein-Milman theorem, $\rho_{\mathcal V}$ is a convex combination of extremal invariant states, so there exists an invariant state $\eta$ such that $\mathrm{supp}\,\eta\subset \mathrm{supp}\, \rho_{\mathcal V}$, and the minimal enclosure $\mathrm{supp}\,\eta$ is contained in $\mathcal V$.

\item If $\mathcal V$ and $\mathcal W$ are minimal enclosures contained in $\mathcal R$, then, as in the proof of point $4$, they are the supports of invariant states $\rho_{\mathcal V}$ and $\rho_{\mathcal W}$. Because the decomposition of $\mathcal V\oplus \mathcal W$ into minimal enclosures is unique, $\rho_{\mathcal V}$ and~$\rho_{\mathcal W}$ are the unique extremal invariant states of $\M_{|\mathcal V \oplus \mathcal W}$. Since the set of invariant states is convex, then by the Krein-Milman theorem, $\rho$ is a convex combination of $\rho_{\mathcal V}$ and $\rho_{\mathcal W}$, so $\rho_{\mathcal C}$ and $\rho_{\mathcal C}'$ must be zero. \fin
\end{enumerate}

\smallskip
We can now return to the study of enclosures generated by vectors of the form~$x\otimes |i\rangle$. Remark that ``non-connectedness of $i$ and $j$ through $\connect$" (Definition~\ref{defi_irroqrw}), when stated in terms of enclosures, is related to the existence of~$x\in\mathfrak h_i$, $y\in\mathfrak h_j$, such that one of the following holds:
\begin{itemize}
\item[\textbf{(a1)}] $ y \otimes \vec j \not\in \enc(x\otimes \vec i)^\perp $ and $ x\otimes \vec i \in \enc (y \otimes \vec j)^\perp$,
\item[\textbf{(a2)}]  $ y \otimes \vec j \in \enc(x\otimes \vec i)^\perp $ and $ x\otimes \vec i \in \enc (y \otimes \vec j)^\perp$.
\end{itemize}

Our first task will be to show that, when restricting to the subspace $\mathcal R$, the situation \textbf{(a1)} cannot appear. The following Lemma indeed holds:

\begin{lemme}\label{coro_conda}
If $x\otimes\vec i$ and $y\otimes \vec j$ are in $\mathcal R$, then one of the following situations holds:
\begin{enumerate}
\item $ x\otimes\vec i\not\in\enc(y\otimes\vec j)^\perp$ and $ y\otimes\vec j\not\in\enc(x\otimes\vec i)^\perp$
\item $\enc(x\otimes\vec i )\perp \enc(y\otimes \vec j)$.
\end{enumerate}
\end{lemme}

\pre

It is sufficient to notice that, if $y\otimes\vec j \in \enc(x\otimes\vec i)^\perp\cap \mathcal R$, then the minimal enclosures containing $x\otimes\vec i$ and $y\otimes\vec j$ are orthogonal.
Indeed, by point~$2$ in Proposition \ref{prop_coherence}, the subspace
$ \enc(x\otimes \vec i)^\perp \cap \mathcal R$ is an enclosure, and it contains~$y\otimes \vec j $ by assumption. 
\fin
\begin{remark}
Beware that, in situation 1 of Lemma \ref{coro_conda}, one may still have $\enc( x\otimes\vec i)$ and $\enc(y\otimes \vec j)$ non-orthogonal but in direct sum, as the following example shows.
\end{remark}

\begin{example}\label{ex_3}
We consider an OQRW $\mathfrak M$ with two sites, \textit{i.e.} $V=\{1,2\}$, and~$\mathfrak h_1=\mathfrak h_2={\mathbb C}^2$, and, for a fixed $p\in]0,1[$,
$$
L_{11}=L_{22}=\sqrt{p\,}\,\id,
\qquad
L_{12}=L_{21}=\sqrt{1-p\,}\, B\;
\qquad \mbox{ with }\quad 
B=\left(\begin{array}{cc}
0&1\\
1&0
\end{array}\right).
$$
We denote the canonical basis of $\mathfrak h$ by $e_1,e_2$. Then
$$\enc(e_1\otimes\vec 1)=\mathrm{Vect}\{e_1\otimes \vec 1, e_2\otimes \vec 2\}$$ 
$$\enc((e_1+e_2)\otimes\vec 1)=\mathrm{Vect}\{(e_1+e_2)\otimes \vec 1, (e_1+e_2)\otimes \vec 2\}$$
are non-orthogonal but have trivial intersection. 
\end{example}

\smallskip

\begin{lemme}\label{lemma_unicitydec}
Let $\mathcal V = E_1\oplus E_2$, where $E_1$ and $E_2$ are minimal enclosures contained in $\mathcal R$. The decomposition of $\mathcal V$ into a direct sum of minimal enclosures is unique if and only if any enclosure $\mathcal W$ such that $\mathcal W\not\perp E_1$ and $\mathcal W\not\perp E_2$ satisfies $\mathcal W \cap \mathcal V = \{0\}$.  If the latter statement holds, then the two enclosures are orthogonal.
\end{lemme}

\pre

Assume the decomposition of $\mathcal V$ as a direct sum of minimal enclosures is unique. Then $E_1 \perp E_2$, otherwise by Proposition \ref{prop_coherence},  \[\mathcal V \cap E_1^\perp = \mathcal V \cap (\mathcal R \cap E_1^\perp  )\] would be an enclosure that does not contain $E_2$, leading to a different decomposition of $\mathcal V$. Now consider a minimal enclosure $\mathcal W$ with $\mathcal W\not\perp E_1$ and $\mathcal W\not\perp E_2$. This implies $\mathcal W\neq E_1$ so by Lemma \ref{lemma_lemenc}, $\mathcal W \cap E_1=\{0\}$. If $\mathcal W\cap \mathcal V\neq \{0\}$ then it is an enclosure in $\mathcal W$ so by minimality, $\mathcal W\subset \mathcal V$. Then $\mathcal W \oplus E_1$ is a direct sum of minimal enclosures contained in $\mathcal V$, so, by point $2$ in Proposition \ref{prop_coherence}, one can complete this as a decomposition of~$\mathcal V$ into a direct sum of minimal enclosures. This is a contradiction, leading to $\mathcal W\cap \mathcal V= \{0\}$.

Now assume that any enclosure $\mathcal W$ such that $\mathcal W\not\perp E_1$ and $\mathcal W\not\perp E_2$ satisfies $\mathcal W \cap \mathcal V = \{0\}$. 
Taking first $\mathcal W= E_2$, which obviously has a non trivial intersection with $\mathcal V$, we obtain that $E_1\perp E_2$.
Now consider some minimal enclosure $E_3$ contained in $\mathcal V$. Then, by assumption, one has \textit{e.g.} 
$E_3\perp E_1$ and $E_3\not\perp E_2$. But then by point $2$ in Proposition \ref{prop_coherence}, one has 
$E_3\subset E_1^\perp \cap \mathcal V$, which, as proved above, is $E_2$. This proves the uniqueness of the decomposition.
\fin
\smallskip

The following remark shows that Lemma \ref{lemma_unicitydec} is consistent with the unicity of the irreducible decomposition for classical Markov chains:
\begin{remark}\label{remark_classicalorthogonal}
Consider a minimal dilation $\M$ of a classical Markov chain. By Proposition \ref{lemma_enctraj} and Lemma \ref{lemma_lemenc}, any minimal enclosure is of the form $\cc~\otimes~\ell^2(V_i)$ for $V_i\subset V$. Therefore, for such an OQRW $\M$, any distinct minimal enclosures~$\mathcal V$ and $\mathcal W$ are always orthogonal.
\end{remark}

\smallskip

Once again, the following result is proven in \cite{BN} in finite dimension. We extend the proof to infinite dimension.
\begin{coro}\label{coro_partialisom}
Assume that $\mathcal V = E_1\oplus E_2$ where $E_1$ and $E_2$ are minimal enclosures contained in $\mathcal R$, but that the decomposition into a direct sum of minimal enclosures, as in Lemma \ref{lemma_unicitydec}, is non-unique. Then 
\begin{equation}\label{eq_egalitedimensions}
\dim\,E_1= \dim \,E_2.
\end{equation}
If, in addition, $E_1\perp  E_2$, 
then there exists a partial isometry $Q$ from $E_1$ to~$E_2$ satisfying
\begin{equation}\label{eq_partialisometry}
Q^* \,Q = \id_{|E_1}\qquad Q\,Q^* = \id_{|E_2}
\end{equation}
and for any $\rho$ in $\mathcal I_1(\mathcal H)$, for $R=Q+Q^*$, and $P_i=P_{E_i}$, $i=1,2$:
\begin{equation}\label{eq_commpartialisom}
R\,\M(\rho)\, P_i + P_i\, \M(\rho)\, R=\M\big(R\,\rho\, P_i + P_i\,\rho\, R\big).
\end{equation}
\end{coro}

\pre

Assume that there exists a minimal enclosure $\mathcal W$ that is distinct from $E_1$ and non-orthogonal to it. Then by point $2$ of Proposition \ref{prop_coherence}, $E_1\cap \mathcal W^\perp$ is an enclosure contained in $E_1$. By minimality of $E_1$ and non-orthogonality between those two enclosures, $E_1\cap\mathcal W^\perp=\{0\}$. Therefore $\dim E_1 \leq \dim\mathcal W$, and by symmetry one has the equality $\dim E_1 = \dim\mathcal W$.

If $E_1\not\perp E_2$, this yields equality \eqref{eq_egalitedimensions}. Otherwise, the non-unicity of the decomposition implies the existence of minimal enclosures $\widetilde E_1$ and $\widetilde E_2$ such that
\[E_1 \oplus E_2=\widetilde E_1 \oplus \widetilde E_2.\]
and one can assume that \textit{e.g.} $\widetilde E_1$ is distinct from both $E_1$ and $E_2$. Necessarily $\widetilde E_1$ is also non-orthogonal to both $E_1$ and $E_2$, and taking $\mathcal W= \widetilde E_1$ we recover equality~\eqref{eq_egalitedimensions}.

Assume now that $E_1\perp E_2$. By the above discussion there exists a minimal enclosure $\mathcal W$ distinct from $E_1$ and non-orthogonal to~it.
Denote by $P_{1}$, $P_{2}$, $P_{\mathcal W}$ the orthogonal projections on $E_1$, $E_2$, $\mathcal W$ respectively. Define the map  $\mathfrak N$ on~$\mathcal B(\H)$ by
\[\mathfrak N : X\mapsto P_{\mathcal R}\,\M^*(X)\,P_{\mathcal R}.\]
One sees immediately that if $E=E_1$, $E_2$ or $\mathcal W$, then $P_{E}$ is (up to multiplication) the unique invariant of $\mathfrak N_{|\mathcal B(E)}$. Consider the decomposition of $P_{\mathcal W}=\begin{pmatrix}A& B^*\\B&C\end{pmatrix}$ in the splitting $\mathcal V=E_1\oplus E_2$, where necessarily $B\neq 0$. A simple consequence of Proposition \ref{prop_coherence} is that in the same decomposition, 
$\N(P_{\mathcal W})=\begin{pmatrix}\N(A)& \N(B)^*\\ \N(B)&\N(C)\end{pmatrix}$. 
Therefore $A$ is proportional to $P_{1}$ and $C$ to $P_{2}$. Writing relations $P=P^*=P^2$ satisfied by $P_{\mathcal W}$, one sees that $B$ must be proportional to an operator $Q$ satisfying the relations~\eqref{eq_partialisometry}. In addition, fixing that same operator $Q$, for $\theta\in[0,\pi]$, the operator that has the form
\[P_{\theta}=\begin{pmatrix}\cos^2 \theta & \sin \theta\cos\theta\, Q^*\\ \sin \theta\cos\theta\, Q & \sin^2\theta \end{pmatrix} \]
is an orthogonal projection preserved by the map $\N$. So its range is an enclosure and, by point $3$ of Proposition \ref{prop_coherence}, $P_\theta$ will satisfy the relation
\[\M(P_{\theta} \, \rho \, P_{\theta})= P_{\theta} \,\M(\rho)\,P_{\theta}.\]
Differentiating this relation with respect to the $\theta$ variable, we have
\[
\M\big(\frac{\mathrm d P_\theta}{\mathrm d\theta} \, \rho \, P_{\theta}+ P_{\theta} \, \rho \, \frac{\mathrm d P_\theta}{\mathrm d\theta}\big)=\frac{\mathrm d P_\theta}{\mathrm d\theta} \, \M(\rho) \, P_{\theta}+ P_{\theta} \, \M(\rho) \, \frac{\mathrm d P_\theta}{\mathrm d\theta}
\] 
Computing the derivatives
at $\theta=0$ and $\theta=\pi/2$, we obtain relation \eqref{eq_commpartialisom}.
\fin
\smallskip
\begin{coro}\label{coro_partialisom2}
Assume that $\mathcal V = E_1\oplus E_2$ where $E_1$ and $E_2$ are mutually orthogonal minimal enclosures, contained in $\mathcal R$, but that the decomposition into a direct sum of minimal enclosures is non-unique. Denote by  $\rhoinv_i$ the unique invariant state with support in $E_i$, $i=1,2$, and by $Q$ the partial isometry defined in Corollary \ref{coro_partialisom}. Then $\rhoinv_2=Q\, \rhoinv_1 \, Q^*$.

If $\rho$ is an invariant state with support in $\mathcal V$,  write $\rho=\begin{pmatrix}\rho_{1,1} & \rho_{1,2}\\ \rho_{2,1} & \rho_{2,2}\end{pmatrix}$. Then:
\begin{itemize}
\item $\rho_{1,1}$ is proportional to $\rhoinv_1$,
\item $\rho_{2,2}$ is proportional to $\rhoinv_2$,
\item $\rho_{1,2}$ is proportional to $\rhoinv_1\, Q^*=Q^*\rhoinv_2$,
\item $\rho_{2,1}$ is proportional to $\rhoinv_2\, Q=Q\rhoinv_1$.
\end{itemize}
\end{coro}

\pre

The first identity is obtained by applying relation \eqref{eq_commpartialisom} to $\rho=\rhoinv_1$ with $P_1$, then applying it again to the resulting relation, this time with $P_2$.

That each $\rho_{i,j}$ is an invariant is an immediate consequence of Proposition \ref{prop_coherence}. The relation satisfied by $\rho_{1,2}$ and $\rho_{2,1}$ is then obtained by applying relation~\eqref{eq_commpartialisom} to \textit{e.g.} $\rho_{1,2}$, with $P_1$ or $P_2$.
\fin

\medskip

We are now in a position to state the relevant decomposition associated to an open quantum random walk $\mathfrak M$.
\begin{prop}\label{prop_finaldec}
Let $\mathfrak M$ be an OQRW on $\H=\bigoplus_{i\in V}\h _i$. There exists an orthogonal decomposition of $\H$ in the form
\begin{equation}\label{eq_finaldec}
\H = \mathcal D \oplus \bigoplus_{\alpha \in A}\enc({x_{\alpha}}\otimes\vec{i_{\alpha}})\oplus \bigoplus_{\beta \in B}\bigoplus_{\gamma \in C_\beta} \enc({x_{\beta,\gamma}}\otimes\vec{i_{\beta,\gamma}})
\end{equation}
such that the sets $A$, $B$, $C_\beta$ are at most countable, $A$ and $B$ can be empty (but not simultaneously), any $C_\beta$ has cardinality at least two, and:
\begin{itemize}
\item every $\enc({x_{\alpha}}\otimes\vec{i_{\alpha}})$ or $\enc({x_{\beta,\gamma}}\otimes\vec{i_{\beta,\gamma}})$ in this decomposition is a minimal enclosure, and therefore an equivalence class for $\biconnect$,
\item for $\alpha$ in $A$, the only minimal enclosure not orthogonal to $\enc({x_{\alpha}}\otimes\vec{i_{\alpha}})$ is $\enc({x_{\alpha}}\otimes\vec{i_{\alpha}})$ itself,
\item for $\beta$ in $B$ and $\gamma\in C_\beta$, any minimal enclosure that is not orthogonal to~$\enc({x_{\beta,\gamma}}\otimes\vec{i_{\beta,\gamma}})$ is contained in $\bigoplus_{\gamma \in C_\beta} \enc({x_{\beta,\gamma}}\otimes\vec{i_{\beta,\gamma}})$.
\end{itemize}
\end{prop}

\pre

We start with the decomposition $\mathcal H = \mathcal D \oplus \mathcal R$, and proceed to decompose~$\mathcal R$. Consider the set of all minimal enclosures $\enc(x\otimes\vec i)$ with the property that the only minimal enclosure non-orthogonal to $\enc(x\otimes\vec i)$ is $\enc(x\otimes\vec i)$ itself. By separability this set is at most countable. We can label these enclosures $\enc(x_\alpha\otimes \vec{i_\alpha})$, $\alpha \in A$. Let $\mathcal O=\bigoplus_{\alpha\in A}\enc(x_\alpha\otimes\vec {i_\alpha})$. Then $\mathcal O$ is an enclosure, and if $\mathcal R \cap \mathcal O^\perp\neq \{0\}$ then, by point $2$ of Proposition \ref{prop_coherence}, it is also an enclosure and we proceed to decompose it. Consider families of minimal enclosures labeled by a set $C$,
$\{\enc({x_\gamma}\otimes\vec{i_\gamma}), \,\gamma\in C\}$ with the property that any minimal enclosure that is not orthogonal to the space $\bigoplus_{\gamma\in C}\enc({x_\gamma}\otimes\vec{i_\gamma})$ is contained in $\bigoplus_{\gamma\in C}\enc({x_\gamma}\otimes~\vec{i_\gamma})$; by the assumption that $\mathcal R \cap \mathcal O^\perp\neq \{0\}$ this set is not empty. Pick a maximal such family, and index it as 
$\{\enc({x_{1,\gamma}}\otimes\vec{i_{1,\gamma}}), \,\gamma\in C_1\}$. 
By point 2 of Proposition \ref{prop_coherence} and Lemma \ref{lemma_lemenc}, one can assume that the different enclosures in this family are mutually orthogonal. If \[\mathcal R\cap \mathcal O^\perp \cap \big(\bigoplus_{\gamma\in C_1}\enc({x_{1,\gamma}}\otimes\vec{i_{1,\gamma}})\big)^\perp\neq \{0\}\] we can iterate this process.
\fin
\begin{remark}\label{remark_classicalorthogonal2}
By Remark \ref{remark_classicalorthogonal} and Lemma \ref{lemma_unicitydec}, any minimal dilation~$\M$ of a classical Markov chain is simply of the form $\H =\mathcal D \oplus \bigoplus_{\alpha \in A}\enc({x_{\alpha}}\otimes\vec{i_{\alpha}})$.
\end{remark}
\smallskip

We will use this decomposition to characterize the form of stationary states. Before we state our next result, let us give some notation. We fix a decomposition \eqref{eq_finaldec} as considered in Proposition \ref{prop_finaldec}. We define for every $\alpha \in A$ and $(\beta,\gamma)\in B\times C_\beta$ the following orthogonal projections (for $\mathcal V$ a subspace of $\H$, the orthogonal projection on $\mathcal V$ is denoted $P_{\mathcal V}$):
\[P_0 = P_{\mathcal D}\qquad P_\alpha = P_{\enc({x_\alpha}\otimes\vec{i_\alpha})}\qquad P_{\beta,\gamma}=P_{\enc({x_{\beta,\gamma}}\otimes\vec{i_{\beta,\gamma}})}\]
and for a state $\rho$, and indices $i$, $j$ taking the values $0$, $\alpha \in A$ or $(\beta,\gamma)\in B\times C_\beta$
\begin{equation}
\label{eq_decomprho} \rho_i=P_i \, \rho \, P_i\qquad  \rho_{i,j}= P_i \, \rho \, P_j.
\end{equation}
When $\mathcal V$ is a subspace of $\mathcal H$ such that $\mathcal I_1(\mathcal V)$ is stable by $\mathfrak M$, we will talk about the restriction $\mathfrak M_{|\mathcal V}$  of $\mathfrak M$ to $\mathcal V$ (instead of the restriction $\mathfrak M_{|\mathcal I_1(\mathcal V)}$ of~$\mathfrak M$ to $\mathcal I_1(\mathcal V)$).
In addition, for $i$ taking the values $\alpha\in A$ or $(\beta,\gamma)\in B\times C_\beta$, we denote by~$\rhoinv_i$ the unique invariant state of $\mathfrak M_{|\enc({x_{\alpha}}\otimes\vec{i_{\alpha}})}$ or $\mathfrak M_{|\enc({x_{\beta,\gamma}}\otimes\vec{i_{\beta,\gamma}})}$.

\begin{theo}\label{theo_invariantstates}
Let $\rho$ be a $\mathfrak M$-invariant state with $\mathcal H$ separable. With the notation~\eqref{eq_decomprho} we have
\begin{enumerate}
\item $\rho_0=0$,
\item every $\rho_{\alpha}$ is proportional to $\rhoinv_\alpha$, which has support $\enc({x_{\alpha}}\otimes\vec{i_{\alpha}})$, 
\item every $\rho_{(\beta,\gamma)}$ is proportional to $\rhoinv_{(\beta,\gamma)}$, which has support $\enc({x_{\beta,\gamma}}\otimes\vec{i_{\beta,\gamma}})$,
\item for $\gamma\neq \gamma'$ in $C_\beta$, the off-diagonal term $\rho_{((\beta,\gamma),(\beta,\gamma'))}$, which we simply denote by $\rho_{(\beta,\gamma,\gamma')}$, may be non-zero, and is an invariant of $\M$. In addition, there exists a partial isometry $Q_{(\beta,\gamma,\gamma')}$ from $\enc({x_{\beta,\gamma}}\otimes\vec{i_{\beta,\gamma}})$ to $\enc({x_{\beta,\gamma'}}\otimes\vec{i_{\beta,\gamma'}})$ such that:
\begin{itemize}
\item $\rhoinv_{(\beta,\gamma')}=Q_{(\beta,\gamma,\gamma')}\, \rhoinv_{(\beta,\gamma)}\, Q_{(\beta,\gamma,\gamma')}^*$
\item $\rho_{(\beta,\gamma,\gamma')}$ is proportional to $Q^*_{(\beta,\gamma,\gamma')}\,\rhoinv_{(\beta,\gamma')}=\rhoinv_{(\beta,\gamma)}\, Q^*_{(\beta,\gamma,\gamma')}$,
\end{itemize}
\item all other $\rho_{i,j}$ (taking the values $0$, $\alpha \in A$ or $(\beta,\gamma)\in B\times C_\beta$) are zero.
\end{enumerate}
\end{theo}

\pre

This follows from a repeated application of Propositions \ref{prop_coherence} and \ref{prop_finaldec}, and Corollary \ref{coro_partialisom2}.
\fin

\begin{remark}
Our main comment here is that there may exist ``coherences" between minimal blocks, \textit{i.e.} non-zero off-diagonal blocks $\rho_{i,j}$, for $i,j$ corresponding to distinct minimal irreducible blocks. Invariant states are not, contrarily to the classical case, just convex combinations of states invariant for the reduced (irreducible) dynamics. We will observe this in Example \ref{example-nonuniquedec}. Note however, that, according to Remark \ref {remark_classicalorthogonal2}, this cannot happen for minimal dilations of classical Markov chains.
\end{remark}

\begin{remark}
One might have hoped that a relevant decomposition of $\M$ would separate sites, \textit{i.e} that one could decompose $\mathcal R$ into a sum of minimal enclosures  $\bigoplus\enc({x_k}\otimes\vec{i_k})$ with  $\enc({x_k}\otimes\vec{i_k}) \subset \bigoplus_{i\in I_k}\mathfrak h_i$ for disjoint $I_k$. This is not true, as Example \ref{example_nodisconnect} shows.
\end{remark}
\begin{example}\label{example_nodisconnect} Consider again Example \ref{ex_2}. We have a unique decomposition of $\H=\mathfrak h~\otimes~\ell^2(V)$ as a sum of minimal enclosures,
$$ \h \otimes \ell^2(V) = \enc({e_1}\otimes\vec{1})\oplus\enc( {e_2}\otimes\vec{1})$$
even though the two minimal enclosures
$$
 \enc( {e_k}\otimes\vec{1})= \cc\,{e_k}\otimes \ell^2(V),
 \quad  k=1,2, 
 $$
connect all three sites. Note also that, in accordance with Lemma \ref{lemma_unicitydec}, the two unique enclosures are mutually orthogonal.
\end{example}

\begin{remark}
Applying Theorem \ref{theo_invariantstates} and the Frigerio-Verri ergodic theorem (see~\cite{FV}) one can obtain results about the ergodic behaviour of $(\M^n(\rho))_n$, that extend  Proposition~\ref{prop_ergodicconvergence} to the reducible case. This will be done in a forthcoming article. However, in certain cases, the results given in the present article can be enough to describe convergence in reducible OQRWs: see Example \ref{ex_apss}.
\end{remark}

\section{Extensions of open quantum random walks}\label{section_extension}

In this section, we define an extension of open quantum random walks, already mentioned in Remark \ref{remark_extensions}. We consider again a countable set of vertices $V$ and a separable Hilbert space $\H = \bigoplus_{i\in V} \h _i$. An extended open quantum random walk will be a map $\tM : \mathcal I_1(\H) \to \mathcal I_1(\H)$ such that if $\rho=\sum_{i,j\in V}\rho(i,j)\otimes \ketbra ij$ then
\begin{equation}\label{eq_extendedoqrw}
\tM(\rho) = \sum_{i\in V} \Big( \sum_{j\in V} \Phi_{i,j}\big(\rho(j,j)\big)\Big)\otimes \ketbra ii
\end{equation}
where each $\Phi_{i,j}$ is a completely positive map from $\mathcal I_1(\h_j)$ to $\mathcal I_1(\h_i)$ such that, for any $j$ in $V$, 
\begin{equation}\label{eq_stochasticphi}
\sum_{i\in V}\Phi_{i,j}^*(\id_{\h_i})= \id_{\h_j}.
\end{equation}
This defines a transition operation matrix in the sense of Gudder (see \cite{Gudder}). Again this $\tM$ maps $\mathcal I_1(\H)$ to the set $\mathcal I_{\mathcal D}$ of block diagonal trace-class operators (see section \ref{section_OQRWs}). In addition, the Kraus representation associates to each $\Phi_{i,j}$ a countable set $E(j,i)$ and, for every $e\in E(j,i)$, a map $L_e$ from $\h_j$ to $\h_i$ such that $\Phi_{i,j}$ can be written as 
\[\Phi_{i,j}(\rho)=\sum_{e\in E(j,i)} L_e^{\,} \,\rho \, L_e^* \quad \mbox{for any }\rho\in \mathcal I_1(\h_j).\]
We view the operators $L_e$ as associated to the edges of a directed multigraph~$(V,E)$ where $E=\cup_{i,j\in V} E(j,i)$.
Then if we denote by $E(j)=\cup_{i\in V}E(j,i)$ the set of outgoing edges at $j$, the stochasticity condition \eqref{eq_stochasticphi} becomes similar to~\eqref{eq_stochastic}:
\[\forall j\in V\quad \sum_{e\in E(j)} L_e^* L_e^{\,} = \id.\]
This reminds us that the present framework encompasses open quantum random walks as defined in the rest of this article. What's more, it should be noted that the power $\M^n$ of an OQRW $\M$ is not in general an OQRW, but is always an extended OQRW.
All the results of the previous sections can be extended to this more general class of evolutions.

\smallskip

As in section \ref{section_OQRWs}, starting from a state $\rho= \sum_{i\in V}\rho(i)\otimes\ketbra ii$ we can define processes ``without measurement" $(\widetilde Q_n,\frac{\tM^n(\rho,\widetilde Q_n)}{\tr\,\tM^n(\rho,\widetilde Q_n)})_{n\in \nn}$: denote 
\[\tM^n(\rho)=\sum_{i\in V}\tM^n(\rho,i)\otimes\ketbra ii. \]
Then the process ``without measurement" is determined by the variable $\widetilde Q_n$, with law
\[\pp(\widetilde Q_n=i)= \tr \,\tM^n(\rho,i)\]
and the process ``with measurement" $(\widetilde X_n, \widetilde \rho_n)_{n\in \nn}$ by
\[(\widetilde X_0,\widetilde \rho_0)=\big(j,\rho(j)\big)  \mbox{ with probability }\tr\,\rho(j)\]
\[ \pp\Big((\widetilde X_{n+1},\widetilde \rho_{n+1})\!=\!(i,\!\frac{\Phi_{i,j}(\widetilde \rho_n)}{\tr\,\Phi_{i,j}(\widetilde \rho_n)}\!)\Big|(\widetilde X_n,\widetilde \rho_n)\!=\!(j,\widetilde\rho_n)\Big)=\tr \,\Phi_{i,j}(\widetilde\rho_n)\quad \forall i\in V.\]
Note that these classical processes associated to $\tM$ were not considered in \cite{Gudder}.
\smallskip

We claim that our vision of open quantum random walks in terms of paths $\pi$ in~$\P(i,j)$ on a directed graph extends to this framework, with paths $\tilde\pi$ in~$\tP(i,j)$ on a directed multigraph.
\smallskip

In particular, we recover all results from sections \ref{section_irreducibility} through \ref{section_stationary}, replacing $\P$ with $\tP$ in our assumptions, and $Q_n, \M^n(\rho,i), X_n, \rho_n$ with $\widetilde Q_n, \tM^n(\rho,i), \widetilde X_n, \widetilde\rho_n$.
More precisely, Proposition \ref{prop_ergodicityOQRW} and Definition \ref{defi_irroqrw} on irreducibility, as well as Lemma \ref{prop_caracaperiodicite} and Theorem \ref{theo_caracaperiodicite} on the period, extend to $\tM$ by simply replacing every $\P$ with $\tP$. Proposition \ref{prop_eqperiodicity} holds if \eqref{prop_eqperiodicity} becomes
\[
P_{k,i} L_e = L_e P_{k\submod 1,j}\quad \forall\, e\in E(j,i).
\]
And similarly Corollary \ref{coro_perturbation} holds if relation \eqref{eq_coroperturbation} becomes
\[  \forall x\in \mathfrak h_i,\ \exists \, e\in E(i,i) \mbox{ such that }\langle x,L_e x\rangle\neq0.
\]
Then, the whole of section \ref{section_recurrence} holds if the processes $Q_n, \M^n(\rho,i), X_n, \rho_n$ are replaced with $\widetilde Q_n, \tM^n(\rho,i), \widetilde X_n, \widetilde\rho_n$. Similarly, sections \ref{section_nonirreducible} and \ref{section_stationary} remain the same, replacing $\P$ with $\tP$ in the definition of enclosures.

\section{Examples and applications}\label{section_examples}

\def\plusn{\overset{\scriptscriptstyle n}{+}}
\def\moinsn{\overset{\scriptscriptstyle n}{-}}

\begin{example}\label{ex_homogeneous} We start with an application to space-homogeneous open quantum random walks on a graph associated with a set of generators of a group. This applies in particular to open quantum random walks on $\zz^d$, which we study in \cite{CP2}.

To be more precise, we assume that $V$ is a set of vertices in an additive (abelian) group $G$, that $\mathfrak h_i = \mathfrak h$ does not depend on $i$, and that there is a finite set $S\subset G$ such that $L_{i,j}=L_{j-i}$ depends only on $j-i$, and is zero unless~$j-i \in S$.
\def\etainv{\eta^{\mathrm{inv}}}

We associate to this OQRW the map 
\begin{equation}\label{eq_defL}
\begin{array}{cccc} \mathfrak L :& \mathcal I (\mathfrak h) & \to & \mathcal I (\mathfrak h) \\ &\eta & \mapsto & \sum_{s\in S} L_s\, \eta \, L_s^* \end{array}.
\end{equation}
If $\M$ is irreducible, then clearly $\mathfrak L$ is also irreducible, and by Proposition \ref{prop_ergodicconvergence}, it has at most one invariant state which we then denote by $\etainv$. Note that, if $\mathfrak h$ is finite-dimensional, then this $\etainv$ exists.
\begin{remark}
From Lemma \ref{lemma_ergodicityCPTPmaps}, one easily sees that $\L$ is irreducible if and only if the operators $L_s$, $s\in S$, have no non-trivial common invariant subspace. This criterion is stated, in particular, in \cite{Fare}.  
\end{remark}

\begin{prop}\label{prop_inexistenceetatinv}
Assume $\M$ as above is irreducible.
\begin{itemize}
\item If $V$ is infinite, then $\M$ does not have an invariant state.
\item If $V$ is finite, then $\mathfrak L$ has an invariant state $\etainv$ and the unique invariant state of $\M$ is  
$$ \sum_{i\in V} \frac{\etainv}{\mathrm{card}\,V} \otimes |i\rangle \langle i |.$$ \end{itemize}
\end{prop}

\pre

Assume there exists an invariant state $\rhoinv$. Since $\M$ is invariant by translation, any translation of that state is also an invariant state, so by Theorem~\ref{theo_unicity}, the state $\rhoinv$ is translation-invariant.
It must therefore be of the form
$\sum_{v\in V} \rho\otimes | v\rangle \langle v|.$
If $V$ is infinite, this has trace either infinite or null and in either case this is a contradiction. If $V$ is finite then it is easy to see that $\rho$ must be an invariant of $\mathfrak L$.
\fin
\smallskip

The Perron-Frobenius theorem for CP maps, Proposition \ref{prop_ergodicconvergence}, allows us to obtain a large deviation principle and a central limit theorem for the position process $(X_n)_{n\in\nn}$ associated with an open quantum random walk $\M$ and an initial state $\rho$ (see section \ref{section_OQRWs}), therefore extending the results of \cite{AGS}. In addition, we can also make more precise the convergence of the sequence of states $(\rho_n)_{n\in\nn}$ (still using the notations of section \ref{section_OQRWs}). This will be done in a separate paper \cite{CP2} studying in detail OQRWs on~$\zz^d$.
\end{example}
\medskip

\begin{example}\label{ex_apss}
We consider the example given in section 12.1 of \cite{APSS}. In our notation this example is given by $V=\{1,2\}$, $\mathfrak h=\cc^2$ (with canonical basis $(e_1,e_2)$) and transitions given by
\[ L_{1,1}=\begin{pmatrix}a & 0 \\ 0 & b\end{pmatrix}\quad  
L_{1,2}=\begin{pmatrix}0 & \!\!\sqrt {p\,} \\ 0 & 0\end{pmatrix}
\quad  L_{2,2}=\begin{pmatrix}1 & 0 \\ 0 & \!\!\sqrt {q\,}\end{pmatrix}\quad  L_{2,1}=\begin{pmatrix}c & 0 \\ 0 & d\end{pmatrix}\]
where we assume $q=1-p\in(0,1)$, $|a|^2+|b|^2=|c|^2+|d|^2=1$, $0<|a|^2,|c|^2<1$. Note that we do not need the additional assumptions $a\neq b$, $c\neq d$, $ab\neq \sqrt{q\,}$, $a^2\neq q$, $b^2\neq q$ done in \cite{APSS}.  
First observe that the only minimal enclosure is 
\[
{\rm Enc}(e_1\otimes \vec 2)
={\rm Vect}(e_1\otimes \vec 2).
\]
Indeed,
\begin{itemize}
\item $\enc(e_1\otimes \vec 1)$ obviously contains 
$\enc(L_{2,1}e_1\otimes \vec 2)=\enc(e_1\otimes \vec 2)$;
\item $\enc(x\otimes \vec 2)$ contains $\enc(L_{1,2}x\otimes \vec 1)$ and if $x=x_1e_1+x_2e_2$ with $x_2\neq 0$, this contains  
$\enc(e_1\otimes \vec 1)$.
\item $\enc(x\otimes \vec 1)$ contains $\enc(L_{2,1}x\otimes \vec 2)
={\rm Enc}\big((cx_1e_1+dx_2e_2)\otimes \vec 2\big)$,
and if $x_2$ is non null, then we fall in the previous case and conclude.
\end{itemize}

Therefore the decomposition \eqref{eq_finaldec} is given by 
\[\mathfrak h \otimes \ell^2(V)= \mathcal D \oplus \big\{\!\begin{pmatrix}a\\0\end{pmatrix}\otimes \vec 2,\ a\in\cc\big\}.\]
By the equivalent definition of $\mathcal D$ given in Remark \ref{remark_RDfinite}, any eigenvector associated to an eigenvalue of modulus one must be orthogonal to $\mathcal D$. So the OQRW $\M$ has a unique eigenvalue of maximum modulus, which is the simple eigenvalue~$1$ associated with the eigenvector
\[\rhoinv=\begin{pmatrix}1 & 0 \\ 0 & 0\end{pmatrix}\otimes \ketbra 22\]
and this implies that, for any initial state $\rho$, one has $\M^n(\rho)\rightarrow\rhoinv$ as $n\to\infty$.
\end{example}

\begin{example}\label{ex_Vn} We consider a family of examples which extends the main example given in \cite{APSS}. This family is indexed by $n\in\nn^*\cup\{\infty\}$; every $\mathfrak h$ is $\cc^2$ and $V$ is either $V_n=\{1,\ldots,n\}$ or $V_\infty=\zz$, and the operators $L_{i,j}$ are defined by
$$ L_{i\addmodn 1,i}=L_+=\frac1{\sqrt 3} \begin{pmatrix} \hphantom{,}1 & 1 \\ \hphantom{,}0 & 1 \end{pmatrix},\qquad 
L_{i\submodn 1,i}=L_-=\frac1{\sqrt 3}\begin{pmatrix} 1 & 0 \\ -1 & 1 \end{pmatrix},$$
where here $\addmodn$, $\submodn$ denote addition or substraction modulo $n$ in the case where $n<\infty$, and standard addition or substraction if $n=\infty$. We denote by $\M_{(n)}$ the above open quantum random walk.

We first show that, in any case, this chain is irreducible, using the characterization given in Proposition \ref{prop_ergodicityOQRW}.
 For this, fix $i$ and $j$ in $\nn$, and let $\Delta=i-k$. For~$p$ large enough, consider $\pi$ of the form $(i,i-1,\ldots,i-\Delta-p,i-\Delta-p+1,\ldots,j)$ (\textit{i.e.} one first moves down $p+\Delta$ times, then up $p$ times), one has 
\begin{eqnarray*}
L_\pi&=&L_+^{p}L_-^{\Delta+p}\\
&=& 3^{-p-\Delta/2} [ \left( \begin{array}{cc} 1 & 0 \\ -\Delta & 1 \end{array}\right) + p \left( \begin{array}{cc} -\Delta & 1 \\ -1 & 0 \end{array}\right) - p^2 \left( \begin{array}{cc} 1 & 0 \\ 0 & 0 \end{array}\right)].
\end{eqnarray*}
Assume that some vectors $x_i=\begin{pmatrix} a_i\\ b_i \end{pmatrix}$ and $x_j=\begin{pmatrix} a_j\\ b_j\end{pmatrix}$ satisfy
$ \langle x_j, L_+^{p}L_-^{\Delta+p}\,x_i\rangle~=~0$
for arbitrarily large $p$. Then one must have
$$ \langle x_j, \begin{pmatrix} 1 & 0 \\ -\Delta & 1 \end{pmatrix}x_i\rangle =  \langle x_j, \left(\begin{array}{cc} \Delta & -1 \\ 1 & 0 \end{array}\right) x_i\rangle =
\langle x_j, \left(\begin{array}{cc} 1 & 0 \\ 0 & 0 \end{array}\right) x_i\rangle =0.$$
By inspection we see that these conditions imply $a_i=b_i=0$ or $a_j=b_j=0$. Therefore, the set $L_\pi x_i$ is total in $\mathfrak h_j$, for any choice of $x_i$.

We now discuss the period. First notice that, for any non null vector $x$ in~$\mathbb C^2$, we always have either 
$\langle x,L_+L_- x\rangle\neq 0$ or $\langle x,L_-L_+ x\rangle\neq 0$. This implies that $D(i,x)\in \{1,2\}$ (just using relation \eqref{def_Dix}) for all $i\in V$ and all $x$, so, by Theorem \ref{theo_caracaperiodicite}, the period can be only $1$ or $2$.

If $n$ is odd, then for $p\in \mathbb N^*$, consider $x=\begin{pmatrix}a\\b\end{pmatrix}\neq 0$. Then
\begin{equation}\label{eq_periodVn} 
\braket{x}{L_+^{pn}\,x}=\braket{\begin{pmatrix}a\\b\end{pmatrix}}{\frac1{3^{pn/2}}\begin{pmatrix}1 & np\\0 & 1 \end{pmatrix}\begin{pmatrix}a\\b\end{pmatrix}} = \frac1{3^{pn/2}} \, (a^2+np\, ab + b^2)
\end{equation}
(this quantity is associated to the path $\pi=1,\ldots, n,\ldots, 1,\ldots, n, 1$ starting from $1$ and going ``up", doing $p$ loops before stopping at $1$). Since $x\neq 0$, the quantity \eqref{eq_periodVn} is zero for at most one $p$, so $D(1,x)$, defined in \eqref{def_Dix}, divides $pn$ for any large enough $p\in\mathbb N^*$. Consequently $D(1,x)=1$ and, by Theorem \ref{theo_caracaperiodicite}, the period is~$1$. By translation-invariance, $D(i,x)=1$ for all $i$ in $V$.

On the other hand, if $n$ is even or infinite, it is clear that the chain has period $2$: the projections
$$ P_{\mathrm{even}}=\sum_{i\,\mathrm{even}} \id\otimes|i\rangle\langle i|\quad \mbox{and} \quad P_{\mathrm{odd}}=\sum_{i\,\mathrm{odd}} \id\otimes|i\rangle\langle i|$$
are $\mathfrak M$-cyclic.

We define one more open quantum random walk, to illustrate the method of ``adding loops" described in Remark \ref{remark_loops} to make an OQRW aperiodic: we define for $\eps\in]0,1[$ the open quantum random walk $\M_{(n,\eps)}$ with sites $V_n$ and transition operators
\[L_{i\addmodn 1,i}^{(\eps)}=L_+^{(\eps)}=\sqrt{1-\eps\,}\,L_+\qquad L_{i\submodn 1,i}^{(\eps)}=L_-^{(\eps)}=\sqrt{1-\eps\,}\,L_-\quad L_{i,i}^{(\eps)}=\sqrt{\eps\,}\,\mathrm{Id}. \]
Note that we consider this perturbation by ``adding a loop" at every site, because it simplifies both the computation of the invariant state, and the simulation. Then $\M_{(n,\eps)}$ is clearly irreducible and, from Corollary \ref{coro_perturbation}, it is aperiodic.

For each choice of open quantum random walk $\M_{(n)}$ (respectively $\M_{(n,\eps)}$) we associate a map $\mathfrak L_{(n)}$ (respectively $\mathfrak L_{(n,\eps)}$) on $\mathcal I_1(\mathbb C^2)$, as in \eqref{eq_defL}. We can check that in all cases, the state $\frac12\,\mathrm{Id}$ on $\mathbb C^2$ is the only invariant of that map. By Proposition \ref{prop_inexistenceetatinv}, for $n\in\mathbb N^*$, the only invariant map of $\M_{(n)}$ (respectively $\M_{(n,\eps)}$) is
\[\rhoinv=\sum_{i\in V_n} \frac1{2n}\,\mathrm{Id}\otimes \ketbra ii \]
We summarize our results:
\begin{prop} Consider the open quantum random walks $\M_{(n)}$ and $\M_{(n,\eps)}$ as above. We have:
\begin{itemize}
\item for every $n$ in $\mathbb N^*\cup\{\infty\}$, the OQRWs $\M_{(n)}$ and $\M_{(n,\eps)}$ are irreducible,
\item for $n$ in $2\,\mathbb N^*\cup\{\infty\}$ the OQRW $\M_{(n)}$ has period 2,
\item for $n$ in $2\,\mathbb N\!+\!1$ the OQRW $\M_{(n)}$ is aperiodic
\item for $n$ in $\mathbb N^*\cup\{\infty\}$, the OQRW $\M_{(n,\eps)}$ is aperiodic,
\item for $n$ in $\mathbb N^*$, the OQRWs  $\M_{(n)}$ and $\M_{(n,\eps)}$ have as unique invariant state
\[\rhoinv=\sum_{i\in V_n} \frac1{2n}\,\mathrm{Id}\otimes \ketbra ii. \]
\end{itemize}
\end{prop}

\smallskip

We now describe the results of numerical simulations. Because we cannot display all data, we choose to focus on what happens ``at site 1". We always start from the initial state
$\rho=\begin{pmatrix}1 & 0 \\ 0 & 0 \end{pmatrix} \otimes |1\rangle \langle 1|$, but the phenomena are insensitive to the particular choice of $\rho$. Whenever we describe a state on $\cc^2$, we give its $(1,1)$ and $(1,2)$ coordinates. Note that:
\begin{itemize}
\item these two coordinates describe the state entirely,
\item because of our choice of $\rho$ and $L_+$, $L_-$, those coordinates are real.
\end{itemize} 
\medskip
In every case, we display for different values of $n$:
\begin{enumerate}
\item the probability $\pp(Q_n=1)$, and its average $\frac1n\sum_{k=0}^{n-1}\pp(Q_k=1)$ (Figures~1,3,5, top row),
\item the $(1,1)$ and $(1,2)$-coefficients of the (non-normalized) ``state at site $1$", \emph{i.e.} $\mathfrak M^n(\rho,1)$ (Figures~1,3,5, middle row), and of the average $\frac1n\sum_{k=0}^{n-1}\mathfrak M^k(\rho,1)$ (Fig. 1,3,5, bottom row),
\item the different values of $X_n$ in a (randomly chosen) quantum trajectory, and the proportion of $1$'s in $X_0,\ldots,X_{n-1}$ (Figures 2,4,6, top row),
\item the $(1,1)$ and $(1,2)$-coefficients of the (normalized) state $\rho_k$ for those times~$k\leq n$ such that $X_k=1$ (Figures 2,4,6, middle row), and of the average $\frac1{N_{n,1}}\sum_{k=0}^{n-1}\rho_k\,\ind_{X_k=1}$ where $N_{n,1}$ is the number of $k$ in $\{0,\ldots,n-1\}$ such that $X_k=1$ (Figures 2,4,6, bottom row).
\end{enumerate}
The series of data 1 and 2 (corresponding to Figures 1,3,5) we call ``without measurement", the series 3 and 4 (corresponding to Figures 2,4,6) we call ``with measurement".
\medskip

\noindent\textbf{Open quantum random walk $\M_{(3)}$}

We obtain numerically the data shown in Figures \ref{V3nonmeasurement} and \ref{V3measurement}. We observe all the convergences mentioned in Corollaries \ref{coro_rec1}, \ref{coro_rec2}, \ref{coro_rec3}.

\begin{figure*}[H]
\hspace{-2cm}\includegraphics[width=1.4\textwidth]{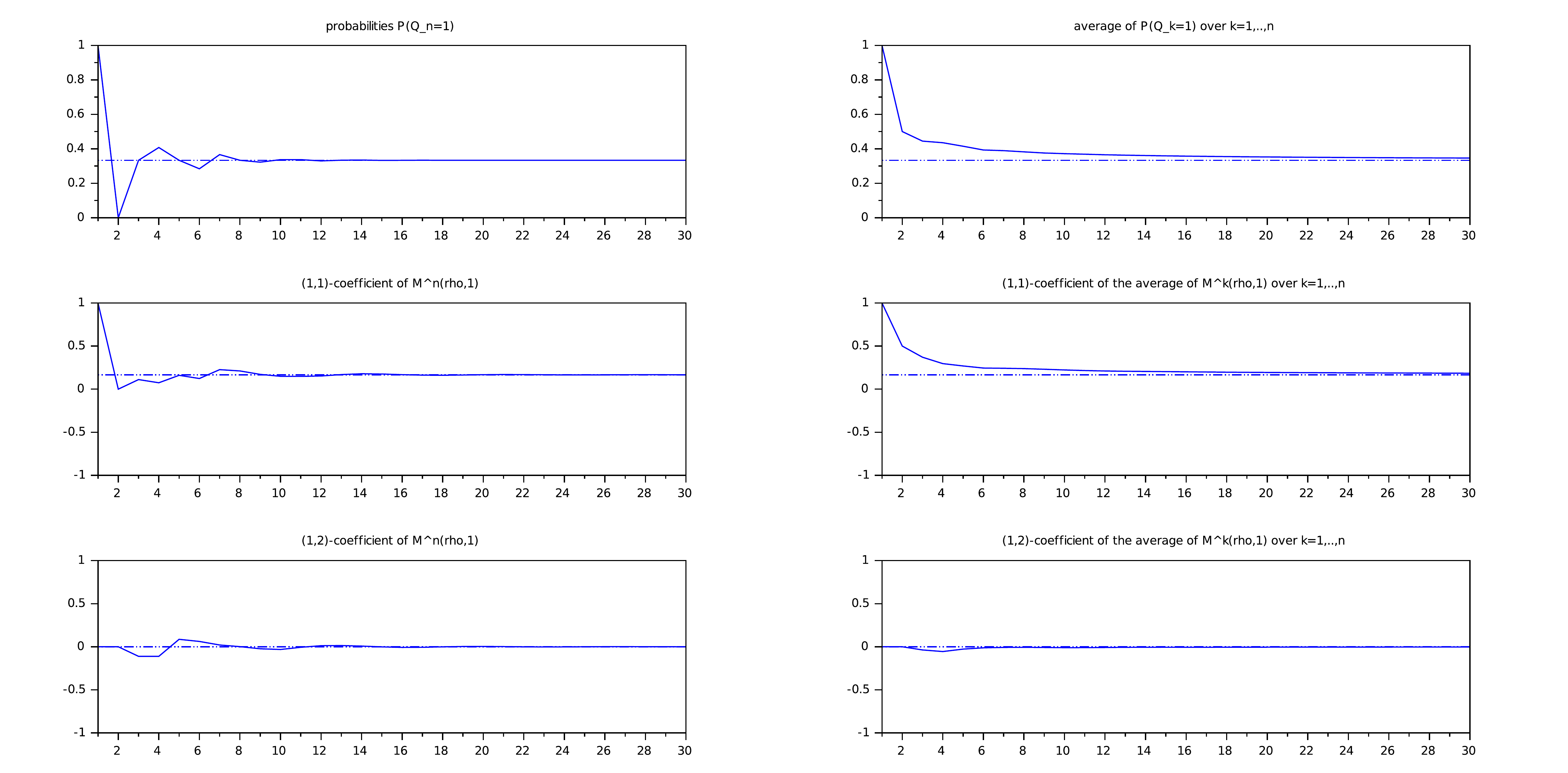}
\caption{OQRW $\M_3$, data without measurements}
\label{V3nonmeasurement}
\end{figure*}

\begin{figure*}[H]
\hspace{-2cm}\includegraphics[width=1.4\textwidth]{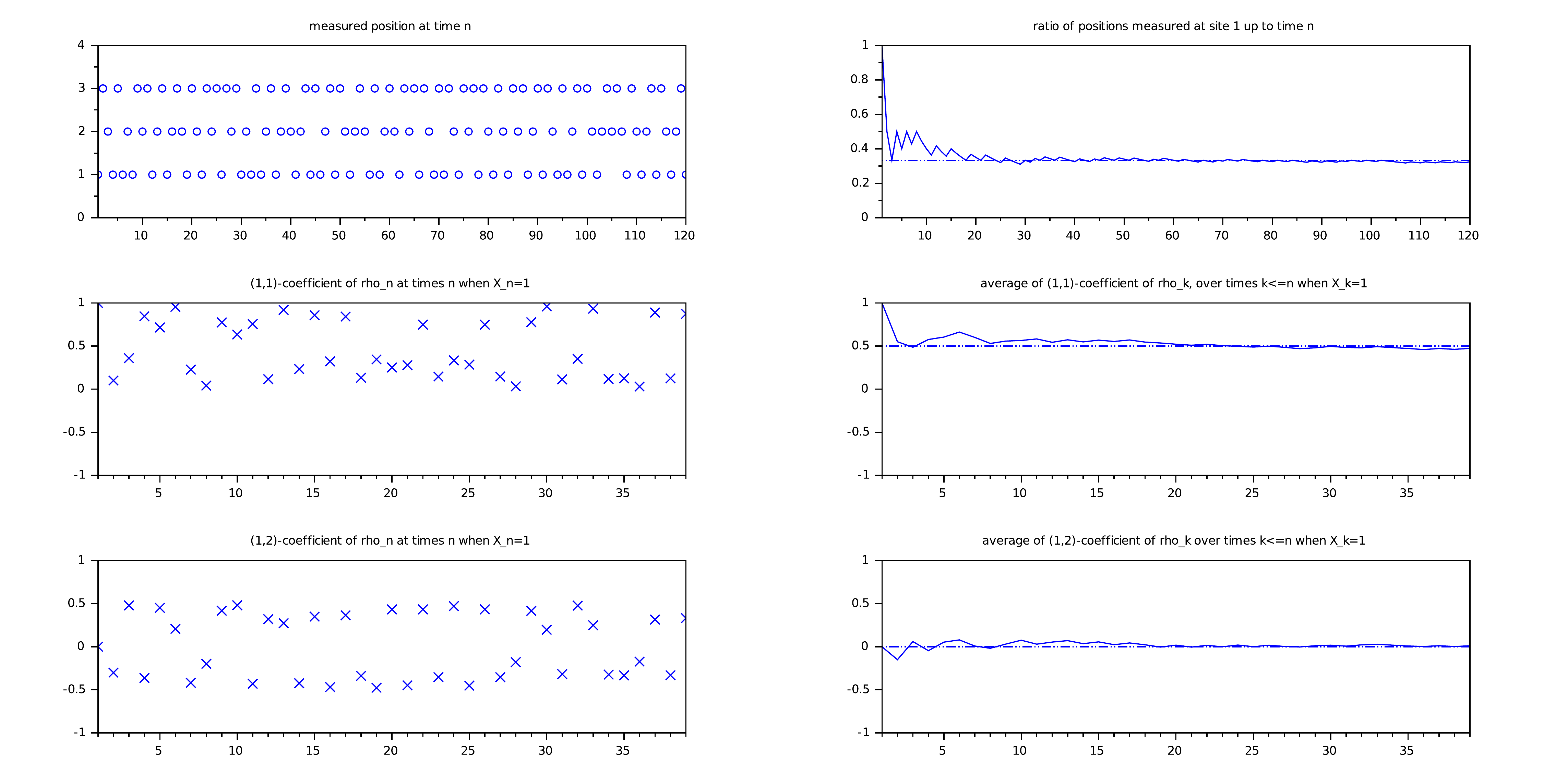}
\caption{OQRW $\M_3$, data with measurements}
\label{V3measurement}
\end{figure*}

\smallskip
\noindent\textbf{Open quantum random walk $\M_{(4)}$}

We obtain numerically the data in Figures \ref{V4nonmeasurement} and \ref{V4measurement}. We observe the convergences mentioned in Corollaries \ref{coro_rec1}, \ref{coro_rec2} but not that of Corollary \ref{coro_rec3}, as the OQRW is not aperiodic. The sequences $\pp(Q_n=1)$ and $\mathfrak M^n(\rho,1)$ exhibit periodic behavior, in a way that is reminiscent of periodic (classical) Markov chains.

\begin{figure*}[H]
\hspace{-2cm}\includegraphics[width=1.4\textwidth]{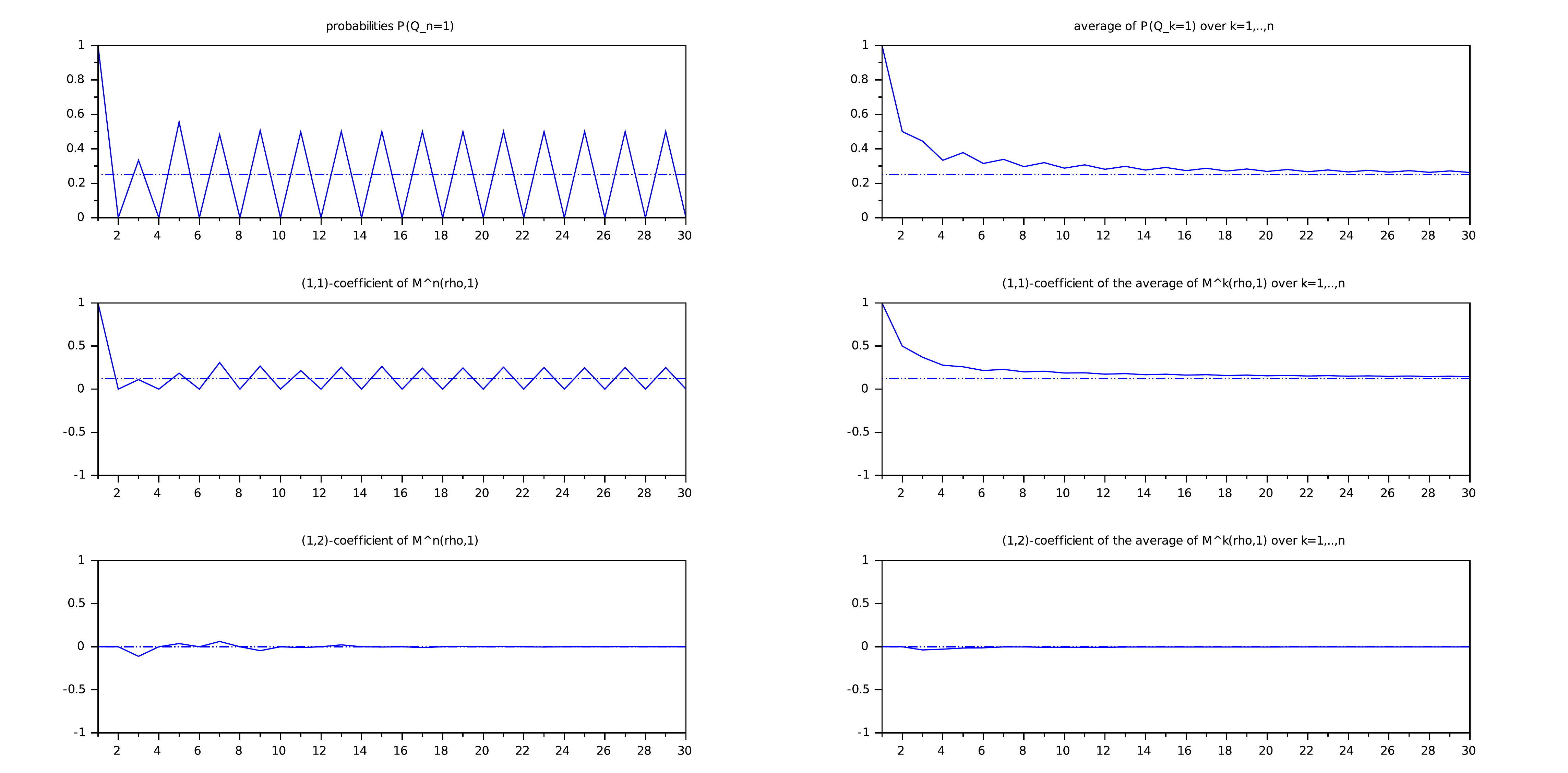}
\caption{OQRW $\M_4$, data without measurements}
\label{V4nonmeasurement}
\end{figure*}

\begin{figure*}[H]
\hspace{-2cm}\includegraphics[width=1.4\textwidth]{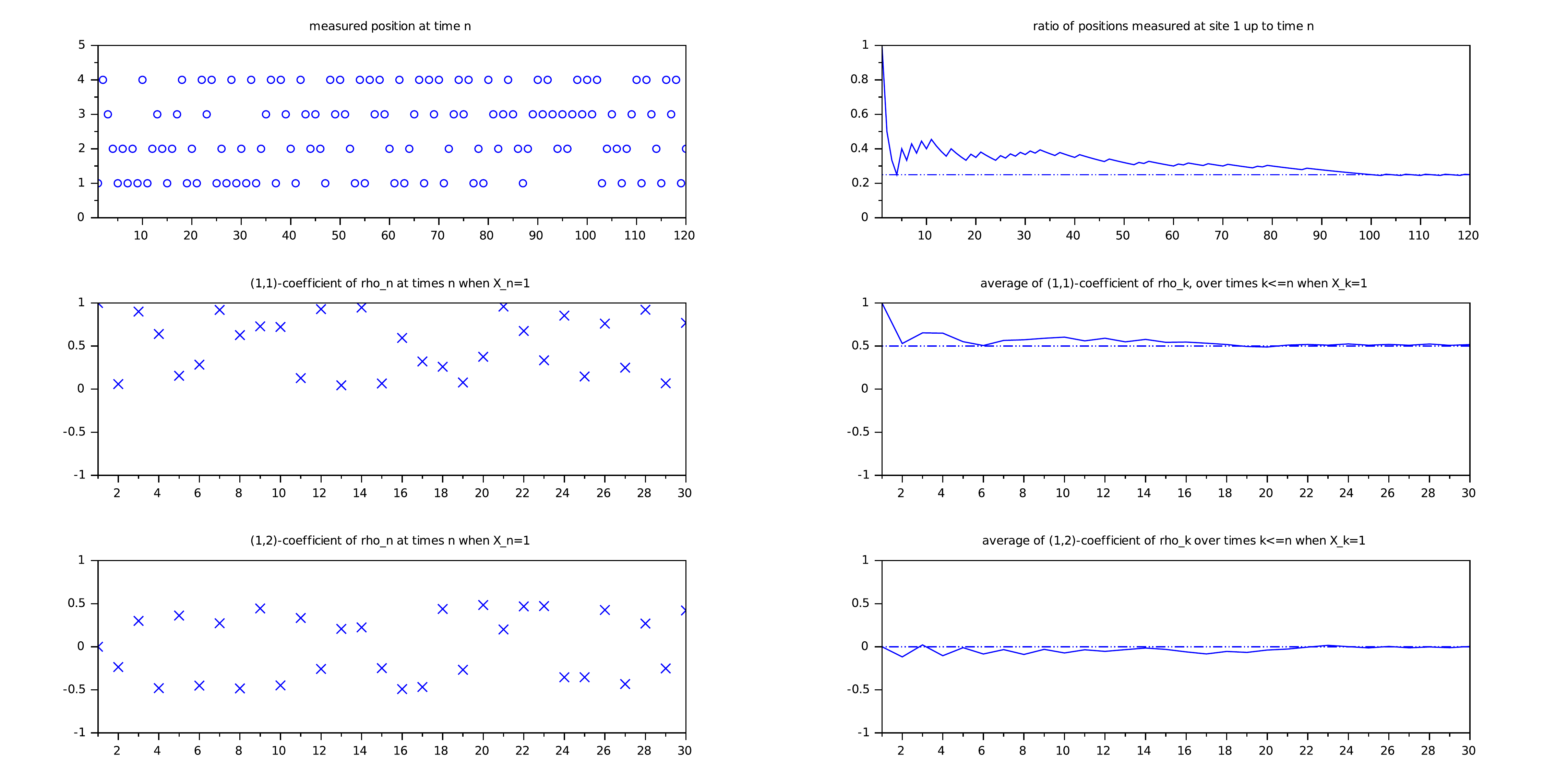}
\caption{OQRW $\M_4$, data with measurements}
\label{V4measurement}
\end{figure*}

\smallskip
\noindent\textbf{Open quantum random walk $\M_{(4,\eps)}$ for $\eps=0.05$}

We obtain numerically the data shown in Figures \ref{V4pertnonmeasurement} and \ref{V4pertmeasurement}. In addition to the convergences mentioned in Corollaries \ref{coro_rec1}, \ref{coro_rec2} we recover those of Corollary~\ref{coro_rec3}, as we have perturbed the OQRW into an aperiodic one.

\begin{figure*}[H]
\hspace{-2cm}\includegraphics[width=1.4\textwidth]{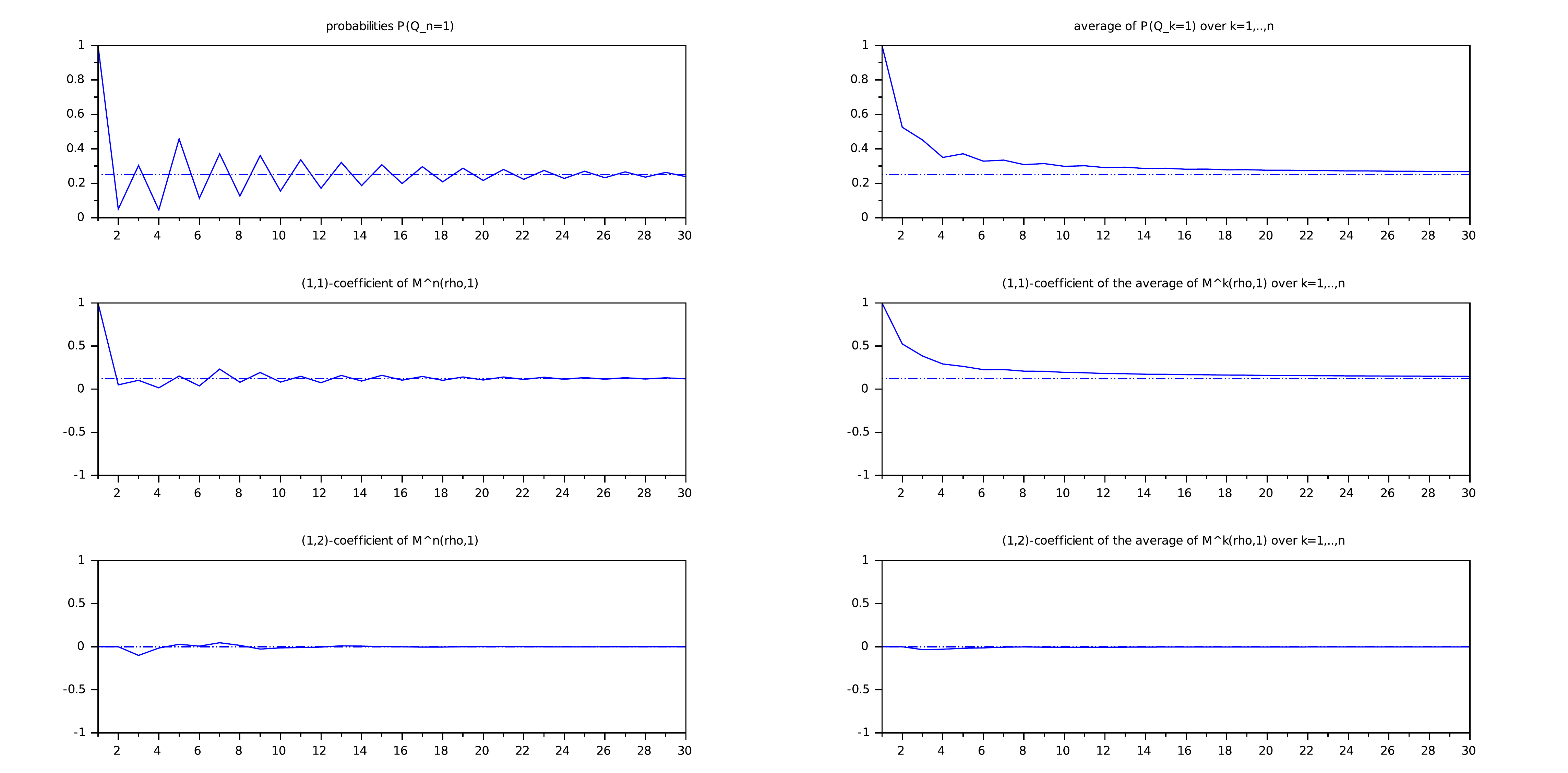}
\caption{perturbed OQRW $\M_{(4,0.05)}$, data without measurements}
\label{V4pertnonmeasurement}
\end{figure*}

\begin{figure*}[H]
\hspace{-2cm}\includegraphics[width=1.4\textwidth]{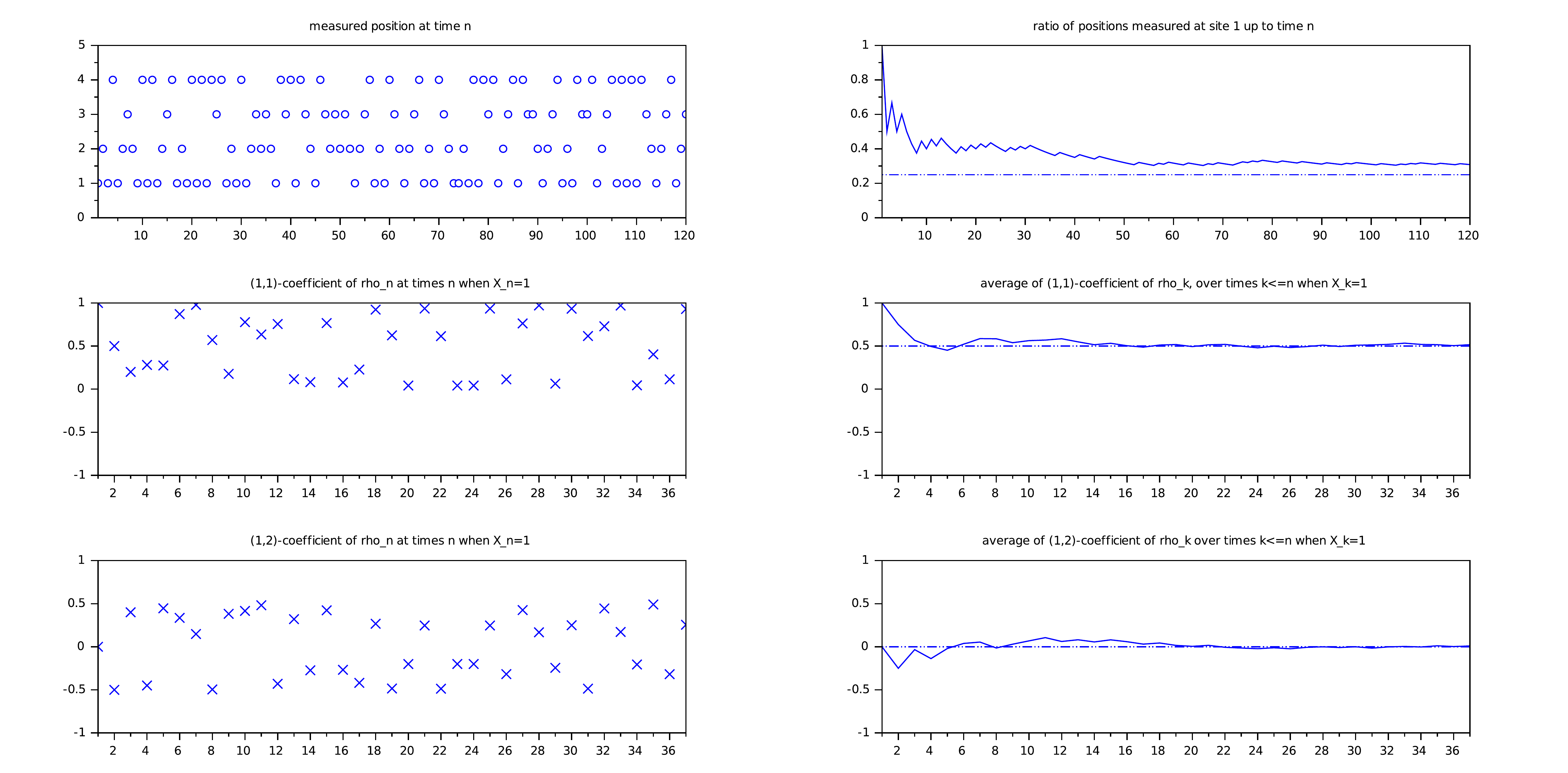}
\caption{perturbed OQRW $\M_{(4,0.05)}$, data with measurements}
\label{V4pertmeasurement}
\end{figure*}

\begin{remark} the data we obtained show that aperiodicity does not imply a convergence of $\rho_n$, even when we condition it on a measurement of $X_n$ at a given site: only convergence in the mean holds. 
\end{remark}
\end{example}

\begin{example}\label{example_Dnotconstant}

We use $V_\infty = \mathbb Z$, $\mathfrak h=\mathbb C^2$ as in the previous example and change the transition matrices,
$$ 
L_+=p\left( \begin{array}{cc} 0 & 1 \\ 1 & 0 \end{array}\right),
\qquad
L_-=q\left( \begin{array}{cc} 1 & 0 \\ 0 & e^{i\alpha} \end{array}\right)
$$
with $\alpha\in[0,2\pi)$, $p,q\in{\mathbb C}\setminus \{0\}$, $|p|^2+|q|^2=1$.

This OQRW is irreducible when $\alpha\neq 0, \pi$. We shall denote by $\{e_0,e_1\}$ the orthonormal basis of $\mathfrak h$ with respect to which we have written the matrix representation of the operators $L_+,L_-$. 
Then it is easy to verify irreducibility by Proposition \ref{prop_ergodicityOQRW} : if we consider the non-zero vector $v=\begin{pmatrix}a\\ b\end{pmatrix}$ in $\mathfrak h_i$, we have that, for all $n>0$, 
\begin{eqnarray*}
&&{\rm span}\{L_+^n v, L_+^{n+1}L_- v, L_+^{n}L_- L_+ v\}
\\
&&=\left\{\begin{array}{ll}
{\rm span}\{\begin{pmatrix}a\\ b\end{pmatrix}, \begin{pmatrix}e^{i\alpha}b\\ a\end{pmatrix} , \begin{pmatrix}b\\ e^{i\alpha} a\end{pmatrix} \} 
& \quad n \, {\rm even},
\\
{\rm span}\{\begin{pmatrix} b\\  a\end{pmatrix} , \begin{pmatrix} e^{i\alpha} a \\ b \end{pmatrix} , 
\begin{pmatrix}a\\ e^{i\alpha} b\end{pmatrix} \} 
& \quad n \, {\rm odd},
\end{array}\right.
\end{eqnarray*} 
 in both cases, it coincides with $\mathfrak h_{i+n}$. Similarly we can proceed for $n\le 0$.  

The period is $4$: we can choose the resolution of the identity
\[P_k = \sum_{i\in \zz} |e_0\rangle \langle e_0|\otimes |4\,i+k\rangle \langle 4\,i+k|
+\sum_{i\in \zz} |e_1\rangle \langle e_1|\otimes |4\,i+k+2\rangle \langle 4\,i+k+2|,\]
for $k=0,\ldots,3$.
Obviously, from the properties of this OQRW and Theorem~\ref{theo_caracaperiodicite}, the period cannot be greater than $4$. So we can conclude that the period is exactly $4$.

Finally, notice that the quantity $D(i,x)$ introduced in Theorem \ref{theo_caracaperiodicite} is not the same for all vectors:
$D(i,e_0)=D(i,e_1)=4$ but, if we call $x=\begin{pmatrix}1\\ e^{i\alpha/2} \end{pmatrix}$, then $x$ is an eigenvector for $L_-L_+$ and so the set of lengths $\ell$ introduced in the definition of~$D(i,x)$ contains $2$. Since it is clear that all those lengths are even, then~$D(i,x)=2$.
\end{example}

\begin{example}\label{example-nonuniquedec}
We consider an OQRW $\mathfrak M$ as introduced in Example \ref{ex_3}. Then~$\mathfrak M$ does not have a unique decomposition in irreducible components. Indeed, it is easy to see that  the $\mathfrak M$-invariant states are all the states of the form
$$
\rho = \rho_1 \otimes |1\rangle\langle 1| + B\rho_1 B \otimes |2\rangle\langle2|
$$
for any $2\times 2$ matrix $\rho_1$ such that $2\rho_1$ is a state in $M_2(\mathbb C)$.
So $\mathcal R=\mathcal H$ for this $\mathfrak M$, and the minimal enclosures are exactly all the enclosures generated by vectors of the form $x\otimes |1\rangle$, for $x =\begin{pmatrix}a\\b\end{pmatrix}$ in $\mathbb C^2$, 
$$
\mbox{Enc}(x\otimes |1\rangle) 
=\mbox{Vect}\{\begin{pmatrix}a\\b\end{pmatrix} \otimes |1\rangle, \begin{pmatrix}b\\a\end{pmatrix} \otimes |2\rangle \}.
$$
Therefore, the decomposition of $\mathcal R =\H$ into a sum of minimal enclosures is non-unique. To illustrate Theorem \ref{theo_invariantstates}, consider an invariant state $\rho$; from the above discussion, it is of the form 
\[\rho = \frac12 \begin{pmatrix}t & s\\ \overline{s}& 1-t\end{pmatrix}\otimes\ketbra 11 +  \frac12 \begin{pmatrix}1-t & \overline{s}\\ {s}& t\end{pmatrix}\otimes\ketbra 22 \]
with $t\in[0,1]$, $|s|^2\leq t(1-t)$. Writing this $\rho$ in the decomposition
\[\H = \enc(\begin{pmatrix}1\\0\end{pmatrix}\otimes \vec 1)\oplus \enc(\begin{pmatrix}0\\1\end{pmatrix}\otimes \vec 1),\]
which is a possible choice of decomposition \eqref{eq_finaldec}, we obtain
\[\rho = \frac12\,\left(\begin{array}{cccc}t & \hphantom{t}0\hphantom{t} & s & 0 \\ \hphantom{t}0\hphantom{t} & t & 0 & s \\ \overline{s} & 0 & 1\!-\!t & 0 \\ 0 & \overline s & 0 & 1\!-\!t\end{array}\right).\]
In agreement with Theorem \ref{theo_invariantstates}, this  $\rho$ is of the form $t\, \rhoinv_1 + (1-t)\, \rhoinv_2 + s\ \eta_{1,2} + \overline{s} \ \eta_{2,1}$, where $\rhoinv_1$ and $\rhoinv_2$ are invariant states with support equal to $\enc(\begin{pmatrix}1\\0\end{pmatrix}\otimes \vec 1)$, $\enc(\begin{pmatrix}0\\1\end{pmatrix}\otimes \vec 1)$ respectively. In addition, the off-diagonal blocks $\eta_{1,2}$ and $\eta_{2,1}$ are also~$\M$-invariant, and with $Q$ the partial isometry of the form 
\[Q=\begin{pmatrix}0&0&0&0 \\ 0&0&0&0 \\ 1&0&0&0 \\ 0&1&0&0\end{pmatrix}\]
we see that $\rhoinv_2=Q\rhoinv_1 Q^* $ and $\eta_{1,2}$ is proportional to $Q^*\rhoinv_2=\rhoinv_1 Q^*$.
\end{example}

\bibliography{biblio}

\end{document}